\documentclass[final,onefignum,onetabnum]{siamart171218}
\pdfoutput=1
\usepackage{amsfonts}
\usepackage{graphicx}
\usepackage{float}
\usepackage{subcaption}
\usepackage{epstopdf}
\usepackage{algorithmic}

\ifpdf
  \DeclareGraphicsExtensions{.eps,.pdf,.png,.jpg}
\else
  \DeclareGraphicsExtensions{.eps}
\fi

\newsiamremark{remark}{Remark}
\newsiamremark{hypothesis}{Hypothesis}
\crefname{hypothesis}{Hypothesis}{Hypotheses}
\newsiamthm{claim}{Claim}

\headers{Gradient Descent-based D-optimal Design}{V. P. Zankin, G. V. Ryzhakov, and I. V. Oseledets}

\title{Gradient Descent-based D-optimal Design for the Least-Squares Polynomial Approximation}

\usepackage{amsopn}

\ifpdf
\hypersetup{
  pdftitle={Gradient Descent-based D-optimal Design for the Least-Squares Polynomial Approximation},
  pdfauthor={V. P. Zankin, G. V. Ryzhakov, and I. V. Oseledets}
}
\fi

\DeclareMathOperator*{\argmax}{\text{arg\hspace*{0.5ex}max}}
\DeclareMathOperator*{\argmin}{\text{arg\hspace*{0.5ex}min}}
\def\vector#1_#2{{#1}_1,\,{#1}_2,\,\ldots,\,{#1}_{#2}}
\def\abs|#1|{\lvert #1\rvert}
\def\norm|#1|{\left\| #1\right\|}
\def\P{\mathsf P}
\def\normi|#1|{\norm|#1|_\infty}
\def\D{{\mathcal D}}
\DeclareMathOperator{\Imm}{Im}

\renewcommand{\vec}[1]{\mathbf{#1}}

\author{Vitaly P. Zankin\footnotemark[1]\ \footnotemark[2]
\and Gleb V. Ryzhakov\footnotemark[1]
\and Ivan V. Oseledets\footnotemark[1]\ \footnotemark[3]}

\begin{document}

\maketitle

\renewcommand{\thefootnote}{\fnsymbol{footnote}}
\footnotetext[1]{Skolkovo Institute of Science and Technology, Skolkovo Innovation Center, Moscow, Russia, (\email{[v.zankin, g.ryzhakov, i.oseledets]@skoltech.ru})}
\footnotetext[2]{Moscow Institute of Physics and Technology, Moscow region, Dolgoprudny, Russia, (\email{zankin@phystech.edu})}
\footnotetext[3]{Institute of Numerical Mathematics of Russian Academy of Sciences, Moscow, Russia}
\renewcommand{\thefootnote}{\arabic{footnote}}

\begin{abstract}
  In this work, we propose a novel sampling method for Design of Experiments. This method allows to sample such input values of the parameters of a computational model for which the constructed surrogate model will have the least possible approximation error. High efficiency of the proposed method is demonstrated by its comparison with other sampling techniques (LHS, Sobol' sequence sampling, and Maxvol sampling) on the problem of least-squares polynomial approximation. Also, numerical experiments for the Lebesgue constant growth for the points sampled by the proposed method are carried out.
\end{abstract}



\section{Introduction}
One of the approaches to the analysis of complex and expensive to evaluate computational models is surrogate modeling. Surrogate modeling methods allow to build a cheap-to-evaluate model that preserves some properties of the initial computational model. 

One of the widely used surrogate modeling methods is Polynomial Chaos Expansion (PCE)~\cite{PolChao}. This method allows to model the response of the original system as a polynomial expansion over some functional basis of orthogonal polynomials. PCE shows impressive results when applied to the system in which inputs are sampled from some probability distribution~\cite{Chkifa2017,Hampton2015,Jakeman2017,Peng2016}. PCE is broadly used as a powerful tool in Uncertainty Quantification~\cite{1703.05312,Guo2018}. In order to build the PCE of some computational model, it is needed to compute the coefficients of PCE. This can be done with the least-squares method. However, the accuracy of the surrogate model that is obtained with the least-squares depends on the so-called~\emph{experimental design} (hereinafter ED). ED is a set of samples~--- points that are taken from the domain of computational model of interest according to some rule. So, the problem of proper selection of ED arises.          

To solve this problem, sampling methods of design of experiments (DoE)~\cite{DoE} are widely used. So-called space-filling designs such as Latin hypercube sampling (LHS) or Sobol' sequence sampling are extensively applied. 
One of the classes of DoE methods is the class of optimal design methods~\cite{Diaz2018,OptDesign}. The methods from this class make it possible to sample such ED that is optimal with respect to some criterion (\textit{e.g.} A-optimality, C-optimality, D-optimality, S-optimality, \textit{etc})~\cite{ACDopt}. The main advantage of the optimal design sampling over the other DoE sampling methods is that for the construction of an accurate surrogate model much smaller number of runs of the initial computational model is required.

In this paper, we propose a new method for sampling of D-optimal~\cite{OptDesign} ED by direct gradient-based optimization of the objective:
\begin{equation}
    \det(A ^T A)\rightarrow\max,
    \label{eq:objective}
\end{equation}
where $A$ is a model matrix that consists of values of all the basis polynomials evaluated at corresponding ED.

We derive an analytical expression for the gradient of \Cref{eq:objective} and test the proposed sampling method on the ordinary least-squares polynomial approximation of the multivariate function. In the context of this testing, an accuracy comparison with LHS, Sobol' sequence sampling, and Maxvol-based sampling~\cite{Goreinov2010,rect_maxvol} is carried out. Also, Lebesgue constant growth~\cite{2014arXiv1407.3291G,Matthias} is investigated.

\section*{Related work} In the paper~\cite{Simpson2001}, a review of statistical techniques for building an approximation of expensive computational codes is conducted. The methods of interest are the design of experiments, response surface methodology, neural networks, kriging. 
This review describes an application of metamodeling techniques in engineering design and the issues with application of statistical methods in deterministic computer experiments.

The scope of the paper~\cite{Giunta2003} covers applications in computational engineering design studies of special DoE techniques. These techniques are designed for deterministic computer simulations and include Hammersley sequence sampling, LHS, and orthogonal array sampling.  Also, pseudo-Monte Carlo sampling and quasi-Monte Carlo sampling were included into consideration.

In the paper~\cite{Goel2008}, two criteria of experimental design are considered. The first criterion allows to reduce the effect of noise during the surrogate construction while the second criterion helps to reduce bias errors. It is stated that a good sampling method should fulfill both criteria at the same time.    
In this paper, multiple criteria for the assessment of widely-used experimental design methods (such as LHS and D-optimal sampling methods) are used. It is demonstrated that the majority of the sampling methods fulfill only one of the criteria but not the other. 

\section{Problem statement}
\label{sec:problem}
Let us consider a computational model describing a certain system (for example, physical) $f (\vec x),$ where $\vec x = (\vector x_d)^T \in \mathcal{X}\subset\mathbb{R}^d$ is the column-vector of the input variables, $y\in\mathbb{R}$ is the output variable, and $\mathcal{X}$ is the set of admissible vectors $\vec x.$

We consider the model $f(\vec x)$ as a black box: we assume that it can be represented in the form of an expansion over a certain basis of orthogonal polynomials:

\begin{equation*}
f(\vec{x}) = \sum\limits_{j\in\mathbb{N}_0}c_j\Psi_j(\vec{x}),\;\; \vec{x}\in\mathcal{X}
\end{equation*}
where $\Psi_j(\vec x)$ is a multivariate polynomial, $c_j$ are expansion coefficients, and $\mathbb{N}_0\equiv\mathbb{N}\cup\{0\}$ is an extended set of natural numbers.

An element of a $d-$dimensional polynomial basis is defined as the tensor product of univariate polynomials:
$$\Psi_{\vec\alpha}(\vec x) = \prod\limits_{i = 1}^d\psi_{\alpha_i}(x_i), \;\; \vec\alpha\in\mathbb{N}^{d}_0, $$
where $\alpha_i$ is the degree of the univariate polynomial.

By choosing the set of multi-indices $\vec\alpha\in\mathcal{A}\subset\mathbb{N}^{d}_0 $ for some rule, we obtain a polynomial expansion of our model of interest:

\begin{equation}\label{eq:1}
f(\vec{x})\approx\widetilde f(\vec{x}) = \sum\limits_{\vec\alpha\in\mathcal{A}}\tilde{c}_{\vec\alpha}\Psi_{\vec\alpha}(\vec x).
\end{equation}

Our goal is by evaluating the function of interest $f(\vec{x})$ at points from its domain to recover the coefficients of the expansion \cref{eq:1}.

In order to truncate the number of terms in the expansion, we will choose the set $\mathcal{A}$ as follows~\cite{Truncation}:

$$\mathcal{A} = \left\{\vec\alpha\in\mathbb{N}^{d}_0\colon\|\vec\alpha\|_q = \sum\limits_{i = 1}^d\alpha_i^q\leq p^q\right\},$$
where $p$ is the total degree of the polynomial, and $q\in(0,1].$

It is easy to see that the cardinality of the set $|\mathcal{A}| $ is decreasing with decreasing of $q$. This truncation scheme allows to decrease the number of terms in polynomial expansion while keeping the same total degree $p$.

In this paper we consider the case $q = 1,$ \emph{i.\,e.} $\|\vec\alpha\|_1 = \sum\limits_{i = 1}^d\alpha_i\leq p.$ 
The cardinality of the set $|\mathcal{A}|$ in this case (or, what is the same, the number of terms of a polynomial expansion) is $
\begin{pmatrix}d + p\\ p\end{pmatrix}.$
Let us define as an \emph{experimental design} (ED) the following matrix $X:$
\begin{equation*}
X = [\vector \vec x_n]^T\in\mathbb{R}^{n\times d},\;\; n\geq|\mathcal{A}|.
\end{equation*}
The model (Vandermonde-like) matrix $A\in\mathbb{R}^{n\times|\mathcal{A}|} $ is defined as:
\begin{equation}
A_{ij} = \prod\limits_{k = 1}^d\psi_{\vec\alpha_j^{(k)}}(
x_i^{(k)}
),
\label{eq:vandermonde}
\end{equation}
where $\vec\alpha_j^{(k)}$ is the $k$-th component of the multi-index $\vec\alpha_j\in\mathcal{A} = \{\vec\alpha_1,\vec\alpha_2,\ldots,\vec\alpha_{|\mathcal{A}|}\},$ and $x_i^{(k)}$ is the $k$-th component of the $i$-th point $\vec{x}_i$ of input. It is worth noting that all the elements of $\mathcal{A}$ are ordered arbitrary and fixed.

The coefficients in \Cref{eq:1} can be found as a solution to the ordinary least-squares minimization problem:
\begin{equation}\label{eq:LSM}
    \tilde{\vec{c}} = A^+\mathcal{Y} = (A^TA)^{-1}A^T\mathcal{Y},
\end{equation}
where $\tilde{\vec{c}} = (\tilde{c}_{\vec\alpha_1},\tilde{c}_{\vec\alpha_2},\ldots,\tilde{c}_{\vec\alpha_{|\mathcal{A}|}})^T$ is a column-vector of coefficients of polynomial expansion, and  $\mathcal{Y} = (\vector y_n)^T$,
$y_i=f(\vec{x}_i)$ is a column-vector of model responses at ED matrix $X$.

We will call a matrix $\widetilde{X}\in\mathbb{R}^{n\times d}$ an \emph{optimal} ED if the following D-optimality criterion holds for it.

\begin{definition}\label{def:dopt}
$\widetilde{X}\in\mathbb{R}^{n\times d} $ is an D-optimal ED if the following criterion is satisfied:
$$\widetilde{X}:\det B(\widetilde{X}) = \max\limits_{X\in\mathbb{R}^{n\times d}}\det B(X),$$
\end{definition}
where $B(X) = A^T(X)\cdot A(X)$ is a symmetric non-negative definite matrix.

\section{Objective function and its gradient}
\label{sec:GD}

To solve the problem posed in \Cref{sec:problem}, we will use the gradient descent method.
Since we use a D-optimality criterion, it is quite a natural way to optimize the following objective function:
$$W(X) = \det B(X),$$
where matrix $B$ is non-negative definite.
Thus, the problem of finding a D-optimal experimental design $\widetilde{X}$ can be written as follows:
$$\widetilde{X} = \argmax\limits_{X\in\mathbb{R}^{n\times|\mathcal{A}|}} W(X)= \argmax\limits_{X\in\mathbb{R}^{n\times|\mathcal{A}|}} \det B(X) = \argmax\limits_{X\in\mathbb{R}^{n\times|\mathcal{A}|}} \det A(X)^TA(X),$$
for a fixed set $\mathcal{A}$ (\textit{i.e.} for a fixed total degree $p$).
For the problem posed in this way, the standard approach is to replace the original problem with the equivalent one~\cite{log_det}:
\begin{equation}\label{problem}
\widetilde{X} = \argmin\limits_{X\in\mathbb{R}^{n\times|\mathcal{A}|}} \widehat{W}(X) = \argmin\limits_{X\in\mathbb{R}^{n\times|\mathcal{A}|}} \log\det B^{-1}(X) = \argmin\limits_{X\in\mathbb{R}^{n\times|\mathcal{A}|}}( -\log\det B(X)).
\end{equation}

It is worth noting that since the function $\widehat{W}(X)$ is differentiable, we can use gradient descent to find the minimum.

\subsection{Analytical calculation of gradient}

We define $G =\nabla_{X}\widehat{W}(X)\in\mathbb{R}^{n\times d}$ as the gradient matrix of the function $\widehat{W}(X).$ The gradient matrix element is 
$$G_{ij} = \dfrac{\partial\widehat{W}(X)}{\partial \vec{x}_i^{(j)}}.$$

Having defined the gradient matrix in such a way, we can obtain an analytical expression for finding its elements $G_{ij}.$
First of all, let us consider a one-dimensional case ($d = 1$). In such a case, the experimental design $X$ is presented as a column-vector of $n$ one-dimensional points
$X = (\vector x_n)^T\in\mathbb{R}^n$.\\ 
Corresponding gradient matrix $G = \left[\dfrac{\partial\widehat{W}(X)}{\partial x_1},\ \ldots,\ 
\dfrac{\partial\widehat{W}(X)}{\partial x_n}\right]^T\in\mathbb{R}^{n\times{1}}.$
The matrix $A(X)\in\mathbb{R}^{n\times(p+1)}$ will have the following form:
\[A(X) = 
\begin{pmatrix}
    \psi_0(x_1) & \psi_1(x_1) & \hdots & \psi_p(x_1) \\
    \psi_0(x_2) & \psi_1(x_2) & \hdots & \psi_p(x_2) \\
    \vdots & \vdots & \ddots & \vdots \\
    \psi_0(x_n) & \psi_1(x_n) & \hdots & \psi_p(x_n)
\end{pmatrix},\]
where $p$ is the total polynomial degree.
Then the matrix $B(X)\in\mathbb{R}^{(p+1)\times(p+1)}$, $B(X)= A(X)^{T}\cdot A(X)$ is symmetric and is represented as follows:
\[B(X) =
\begin {pmatrix}
    \sum\limits_{m = 1}^{n}\psi_0^2(x_m) & \sum\limits_{m = 1}^{n}\psi_0 (x_m)\psi_1(x_m) & \hdots & \sum\limits_{m = 1}^{n}\psi_0(x_m)\psi_p(x_m) \\
    \sum\limits_{m = 1}^{n}\psi_1(x_m)\psi_0(x_m) & \sum\limits_{m = 1}^{n}\psi_1^2(x_m) & \hdots & \sum\limits_{m = 1}^{n}\psi_1(x_m)\psi_p (x_m) \\
    \vdots & \vdots & \ddots & \vdots \\
    \sum\limits_{m = 1}^{n}\psi_n(x_m)\psi_0(x_m) & \sum\limits_{m = 1}^{n}\psi_n(x_m)\psi_1(x_m) & \hdots & \sum\limits_{m = 1}^{n}\psi_n^2(x_m)
\end{pmatrix},\]
where $ B = \sum\limits_{m = 1}^{n}A_{mi}A_{mj} = \sum\limits_{m = 1}^{n}\psi_i(x_m)\psi_j(x_m).$

Let us formulate the Lemma that allows us to calculate the gradient of the objective.  

\begin{lemma}\label{lemma:one}
For the one-dimensional case $(d=1)$, $k$-th $(k\in\{1,\,2\, \ldots,\, n\})$ component of the gradient matrix $G$ is equal to:
$$ G_k = \dfrac{\partial\widehat{W}(X)}{\partial x_k} = -\sum\limits_{ij}(B^{-1}(X))_{ji}\cdot\left[\dfrac{\partial\psi_i(x_k)}{\partial x_k}\cdot\psi_j(x_k) + \psi_i(x_k)\cdot\dfrac{\partial\psi_j(x_k)}{\partial x_k}\right].$$
\end{lemma}

Now, we generalize the result obtained in \cref{lemma:one} to the multidimensional case when $d > 1$. Let us also recall that $\vec\alpha_j^{(k)}$ is the $k$-th component of the multi-index $\vec\alpha_j\in\mathcal{A}.$

\begin {theorem*}\label{th:one}
The element of the matrix $G$ for $d > 1$ is expressed as follows:
\begin{multline*}
G_{kl}  =  -\sum\limits_{ij}(B^{-1}(X))_{ji}\cdot\left[\dfrac{\partial \psi_{\vec\alpha_i^{(k)}}({x}_k^{(l)})}{\partial {x}_k^{(l)}}\cdot\psi_{\vec\alpha_j^{(k)}}({x}_k^{(l)}) + \psi_{\vec\alpha_i^{(k)}}({x}_k^{(l)})\cdot\dfrac{\partial \psi_{\vec\alpha_j^{(k)}}({x}_k^{(l)})}{\partial {x}_k^{(l)}}\right] \\
\times\prod\limits_{\substack{q=0\\(q\neq{k})}}^{d-1}\psi_{\vec\alpha_i^{(q)}}({x}_k^{(q)})\cdot\psi_{\vec\alpha_j^{(q)}}({x}_k^{(q)}).
\end{multline*}
\end{theorem*}
The proofs of \Cref{lemma:one} and \Cref{th:one} are in \cref{sec:appendix}.

\cref{th:one} allows to calculate the gradient of the objective function $\widehat{W}(X).$ This gradient, in turn, can be used in any of the algorithms of gradient descent.

\subsection{Block-coordinate gradient descent heuristic}
We can compute ED faster by replacing the full gradient with another descent direction. At each iteration of the gradient descent algorithm, we change only the coordinates corresponding to one $d$-dimensional point.

We construct the gradient matrix $G$ at the step $k$ of gradient descent as follows. Using \cref{th:one}, we compute $G(X^{(k)})$, and then we choose a row $l$ of the matrix $G(X^{(k)})$ such that the following condition holds:
\begin{equation}\label{eq:line}
l = \argmax_i\sum\limits_{j = 1}^{d}|G_{ij}(X^{(k)})|.    
\end{equation}
Now we change the matrix $G(X^{(k)})$ by setting the elements of the remaining $n-1$ rows equal to zero:
$$
G(X^{(k)})
:= 
\begin{pmatrix}
    0 & 0 & \hdots & 0\\
    G_{l1}(X^{(k)}) & G_{l2}(X^{(k)}) & \hdots & G_{ld}(X^{(k)})\\
    \vdots   & \vdots   & \ddots & \vdots  \\
    0 & 0 & \hdots & 0
\end{pmatrix}.
$$
The matrix $G(X^{(k)})$ obtained in such a way is the descent direction that we use.
Let us consider the calculation of the gradient matrix element $G$ at the step $k + 1$ of the gradient descent in more detail.
According to \cref{th:one}:
$$G_{ml}(X^{(k+1)}) = -\sum\limits_{ij}(B^{-1}(X^{(k+1)}))_{ji}\cdot\dfrac{\partial B_{ij}(X^{(k+1)})}{\partial X^{(k+1)}_{ml}}.$$
Since the matrices at the $k$-th and at the $(k + 1)$-th steps differ in exactly one row $l$ (see \Cref{eq:line}), we can simplify the calculation of the inverse matrix $B^{-1}(X^{(k+1)})$ with the use of previously computed $B^{-1}(X^{(k)}).$
Inspired by ideas of maxvol~\cite{Goreinov2010}, we can do this by using Sherman-Morrison-Woodbury formula~\cite{SWM}:

\begin{equation}\label{eq:SWM}
    B^{-1}(X^{(k+1)}) = B^{-1}(X^{(k)}) - B^{-1}(X^{(k)})U[I_{2\times 2}+VB^{-1}(X^{(k)})U]^{-1}VB^{-1}(X^{(k)}),    
\end{equation}
where $U\in\mathbb{R}^{|\mathcal{A}|\times 2}$ and $V\in\mathbb{R}^{2\times |\mathcal{A}|}$:
$$U
= 
\begin{pmatrix}
    A_{l1}(X^{(k+1)}) & A_{l1}(X^{(k)}) \\
    A_{l2}(X^{(k+1)}) & A_{l2}(X^{(k)}) \\
    \vdots   & \vdots  \\
    A_{l|\mathcal{A}|}(X^{(k+1)}) & A_{l|\mathcal{A}|}(X^{(k)})
\end{pmatrix},
$$
$$V
= 
\begin{pmatrix}
    A_{l1}(X^{(k+1)}) & A_{l2}(X^{(k+1)}) & \hdots & A_{l|\mathcal{A}|}(X^{(k+1)}) \\
    -A_{l1}(X^{(k)}) & -A_{l2}(X^{(k)})   & \hdots & -A_{l|\mathcal{A}|}(X^{(k)})
\end{pmatrix}.
$$
All the steps described above are summarized in \cref{alg:block_coordinate}.

Thus, for each computation of the gradient matrix $G$ it is necessary to calculate~\cref{eq:SWM} that requires $\mathcal{O}(|\mathcal{A}|^2)$ operations in comparison with the calculation of the full-size matrix $B^{-1}$ that takes $\mathcal{O}(|\mathcal{A}|^3)$ operations.

\begin{algorithm}
    \caption{Block-coordinate gradient descent}
    \label{alg:block_coordinate}
    \begin{algorithmic}[1]
        \STATE{Initialize $X^{(0)}\in\mathcal{D}\subset\mathbb{R}^{n\times|\mathcal{A}|}$}
        \STATE{Calculate $B^{-1}(X^{(0)})$}
        \STATE{Using \cref{th:one}, calculate $G(X^{(0)})$}
        \STATE{Find a row with the largest $L_1$-norm $l := \argmax_i\sum\limits_{j = 1}^{d}|G_{ij}(X^{(0)})|$.}
        \STATE{$G_{ij}(X^{(0)}) := 0$ for $\forall i\neq l,j$}
        \STATE{$X^{(1)}\leftarrow\ $GRADIENT DESCENT($X^{(0)}, G(X^{(0)})$)}
        \STATE{Initialize $k := 1$}
        \WHILE{NOT CONVERGED}
            \STATE{$B^{-1}(X^{(k)}) := B^{-1}(X^{(k-1)})\left(I - U[I_{2\times 2}+VB^{-1}(X^{(k-1)})U]^{-1}VB^{-1}(X^{(k-1)})\right)$}
            \STATE{Using \cref{th:one}, calculate $G(X^{(k)})$}
            \STATE{Find a row with the largest $L_1$-norm $l := \argmax_i\sum\limits_{j = 1}^{d}|G_{ij}(X^{(k)})|$.}
            \STATE{$G_{ij}(X^{(k)}) := 0$ for $\forall i\neq l,j$}
            \STATE{$X^{(k+1)}\leftarrow\ $GRADIENT DESCENT($X^{(k)}, G(X^{(k)})$)}
            \STATE{Update $k := k + 1$}
        \ENDWHILE
        \RETURN $X^{(k)}$
    \end{algorithmic}
\end{algorithm}

\section{Numerical experiments}
\label{sec:numerical}
\subsection{Setting}
In this section, we assess the efficiency of the proposed sampling method (denoted as \textit{GD sampling}) by conducting a comparative study with other sampling methods in terms of approximation accuracy and Lebesque constant growth. 

In \Cref{sec:error}, comparison of the accuracy of the least-squares polynomial approximation build on the sampled points is carried out. The proposed sampling method is tested on four analytical models with varying complexity and input dimensionality. They include three two-dimensional analytical functions namely:
\begin{itemize}
    \item Rosenbrock function,
    \item Sine-cosine product function (denoted as \emph{sincos}),
    \item Gaussian function.
\end{itemize} 
Also, the proposed sampling method is tested on the Piston simulation function that is effectively seven-dimensional. Since all of the testing models are analytical (consequently, cheap to evaluate), the appropriate relative approximation error on a test set is then utilized to evaluate the accuracy of the resulting polynomial expansions.

Despite the fact that the model matrix $A$ can be constructed over different polynomial bases (\textit{e.g.} Legendre polynomials or Hermite polynomials), in the numerical experiments below we are considering Chebyshev polynomials as the basis functions without\- loss of generality:
$$
\begin{aligned}
        \psi_i(x) &= T_i(x),\\ 
        T_{i+1}(x) &= 2x\cdot T_i(x) - T_{i-1}(x),
\end{aligned}
$$
with $ T_1(x) = x$, $T_0(x) = 1$.

Since we consider sampling methods that have stochastic nature, the corresponding EDs that are chosen by these sampling methods are affected. In order to take this fact into the account, we run each analysis 50 times. The aim of repetitions is to assess the effect of stochastic variations.

We consider the following sampling techniques:
\begin{itemize}
    \item LHS~\cite{McKay1979},
    \item Sobol' sequence~\cite{Sobol1967}, 
    \item Maxvol-based~\cite{Goreinov2010,rect_maxvol}.
\end{itemize}

\begin{figure}[ht]
\centering{
    \begin{subfigure}[t]{0.32\linewidth}
        \centering{\includegraphics[width=1\linewidth]{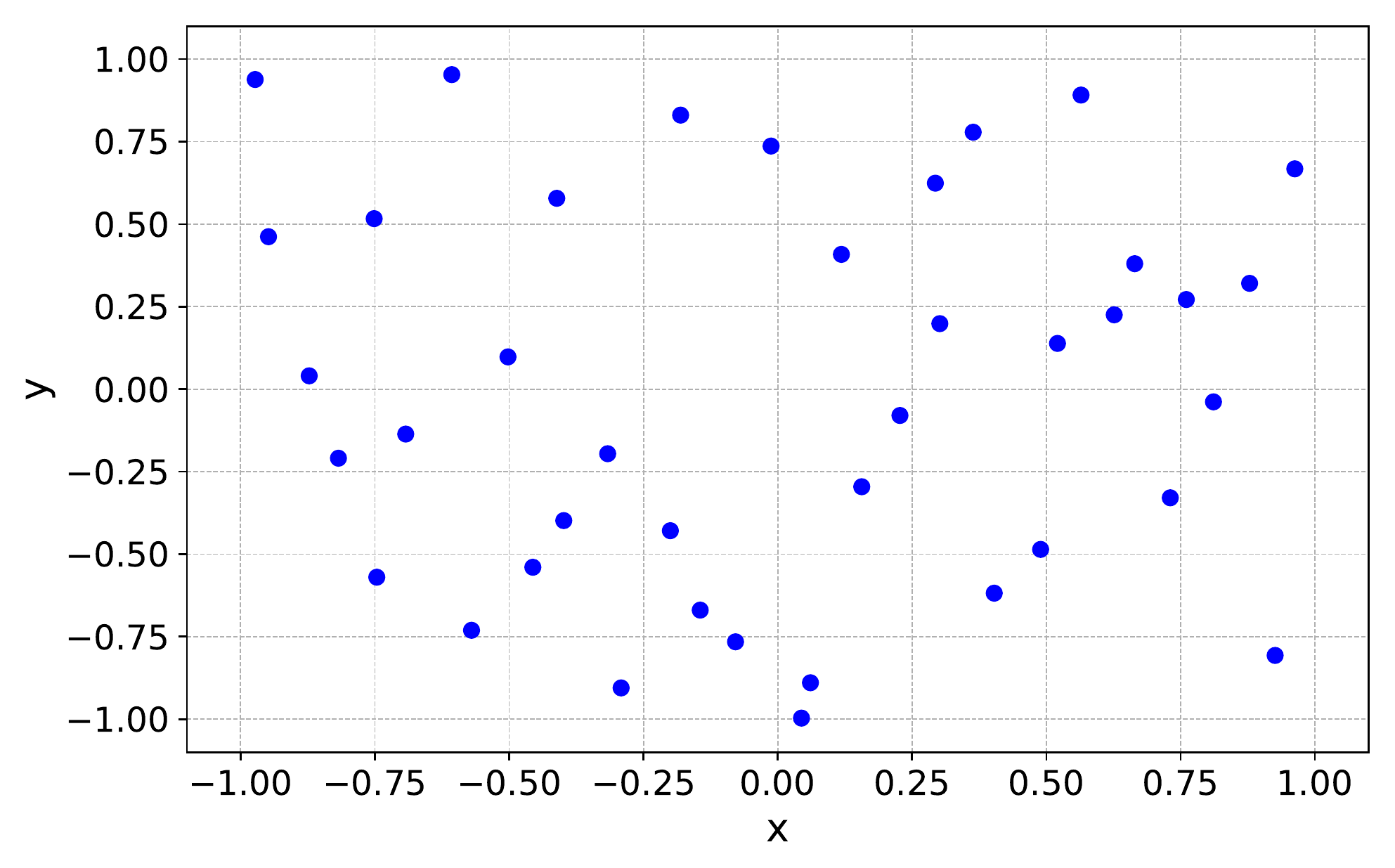}}\caption{LHS}
    \end{subfigure}
        \hfill
    \begin{subfigure}[t]{0.32\linewidth}
        \centering{\includegraphics[width=1\linewidth]{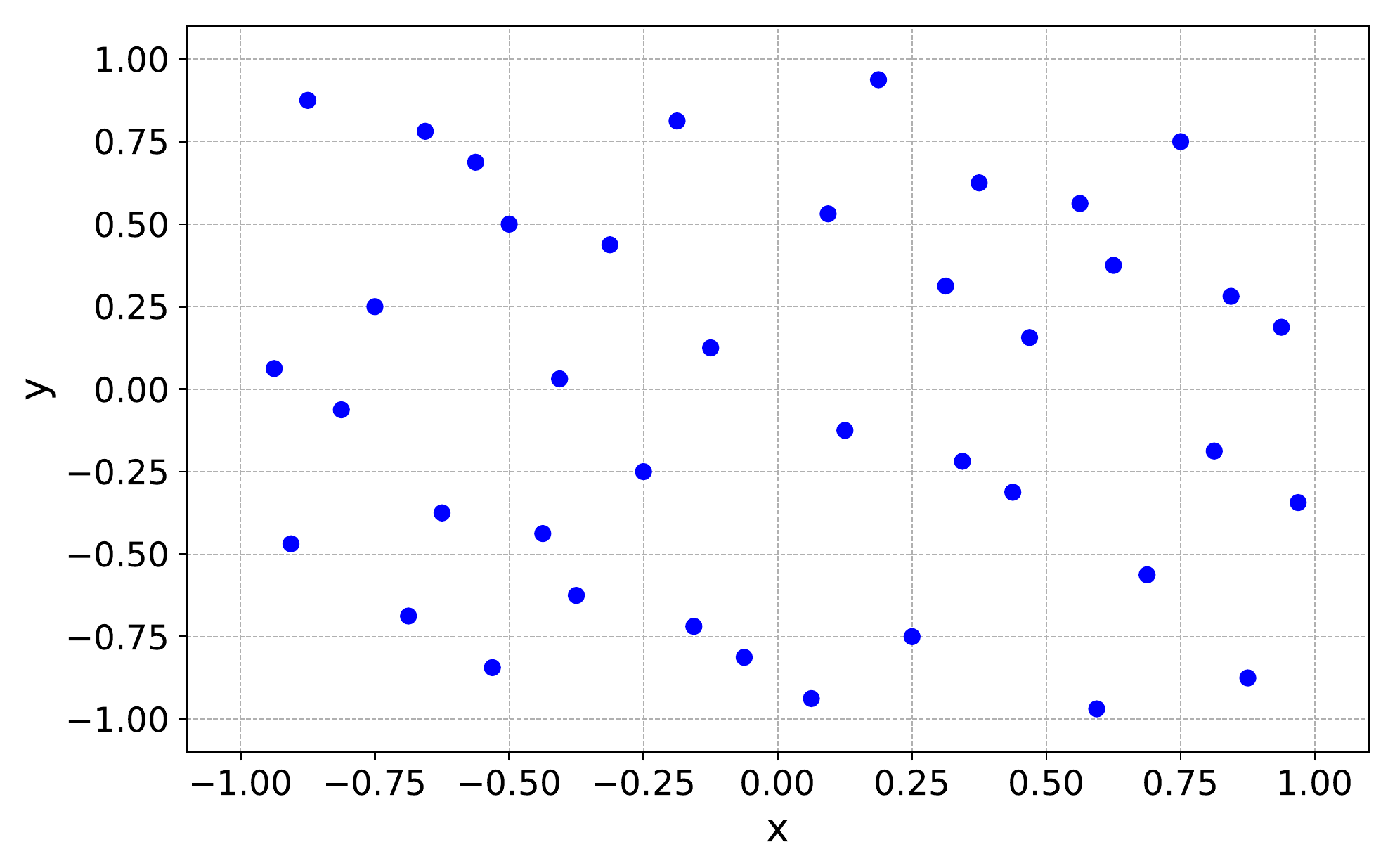}}\caption{Sobol' sequence.}
    \end{subfigure}
    	\hfill
    \begin{subfigure}[t]{0.32\linewidth}
        \centering{\includegraphics[width=1\linewidth]{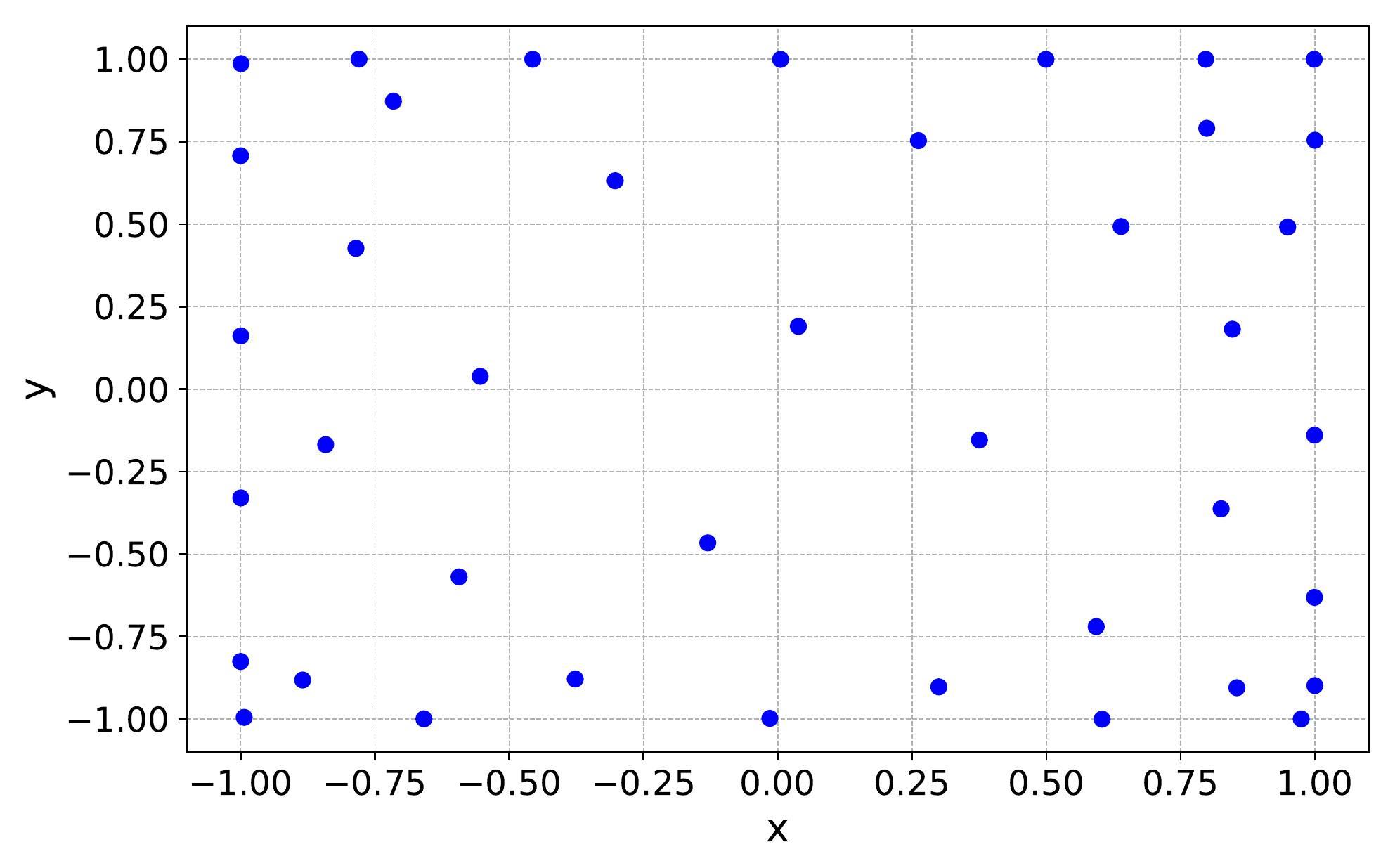}}\caption{Maxvol sampling.}
    \end{subfigure}}\caption{\label{pic:ed_sample} ED of the size 40 obtained with different sampling methods.}
\end{figure}

Examples of ED obtained with the mentioned above sampling techniques for the two-dimensional case are shown at \cref{pic:ed_sample}.

\subsection{Error of approximation}\label{sec:error}
As an accuracy measure of the least-squares polynomial approximation built on ED obtained with each of the sampling techniques, relative error in the infinity norm is utilized:
$$\delta_\infty = \frac{\|f - \tilde{f}\|_\infty}{\|f\|_\infty},$$
where $$\|f\|_\infty\equiv\max\limits_{\vec x\in\mathcal{D}}|f(\vec x)|,$$
for some test set of points  
$\mathcal{D}\subset\mathbb{R}^{d}$.
For all of the experiments below, $\delta_\infty$ is calculated on the test set~$\mathcal D$ of the size $N_{\text{test}}=|\mathcal D|= 10^6$. As a numerical optimizer for GD sampling, we use BFGS method from \texttt{scipy.optimize}.  

\subsubsection{Rosenbrock function}
First of all, we will approximate with a polynomial expansion a well-known two-dimensional Rosenbrock function:
\begin{equation}
f(x,y) = \left(1-x\right)^2 + 100\left(y-x^2\right)^2.
\label{eq:rosenbrock}
\end{equation}

Recall that $l$ is the number of terms in polynomial expansion (in one-dimensional case, it is a total degree minus one of such an expansion), and $n$ is the number of points that make up the ED matrix $X$. In the experiments below, we consider the case $n=l$ when the number of points equals the number of term in polynomial expansion (in such a case, model matrix $A$ is square).

\begin{figure}[htb]
\centering{
    \begin{subfigure}[t]{0.4\linewidth}
        \centering{\includegraphics[width=1\linewidth]{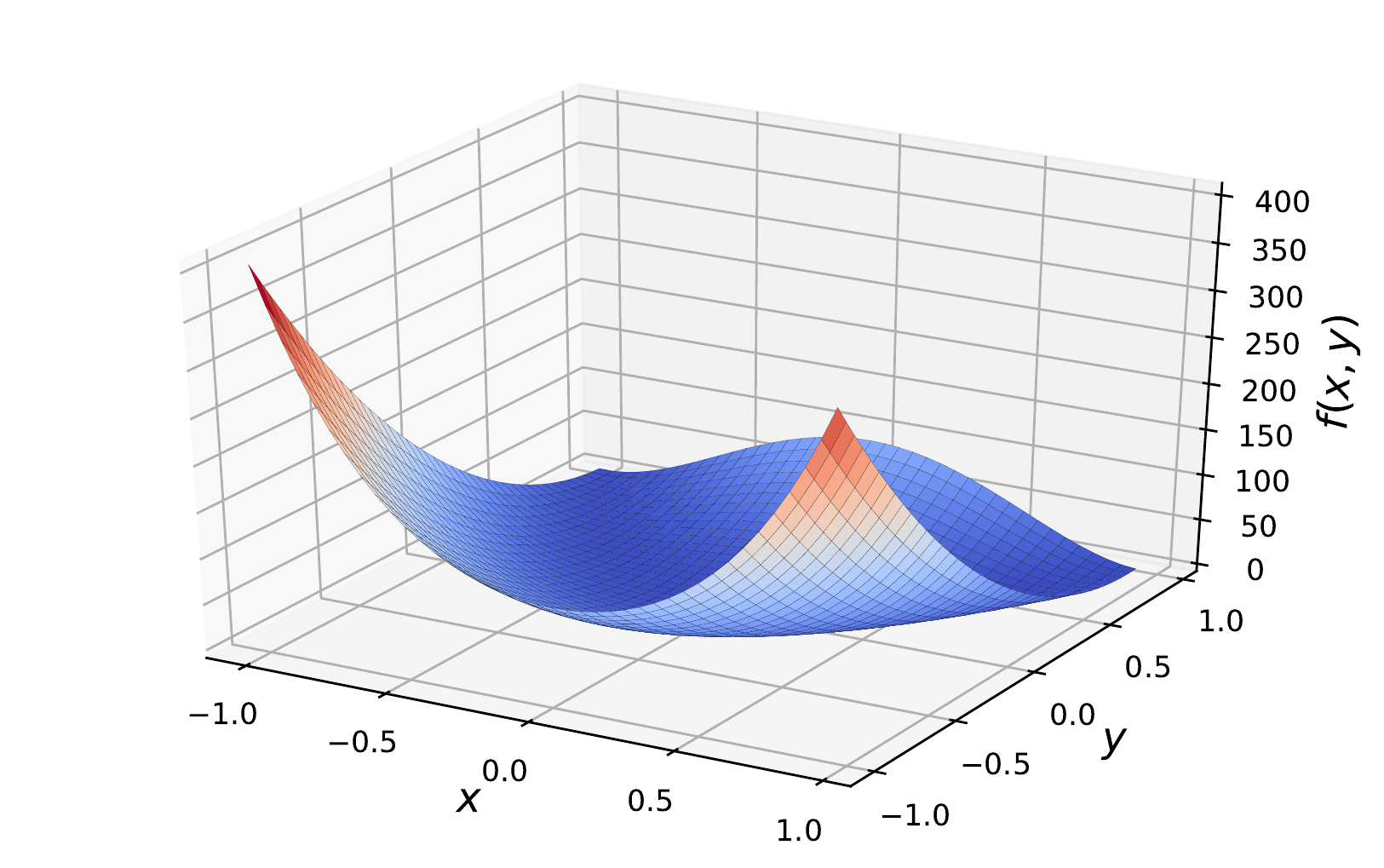}}\caption{Plot of Rosenbrock function.}
    \end{subfigure}
    \vfill
    \begin{subfigure}[t]{0.49\linewidth}
        \centering{\includegraphics[width=1\linewidth]{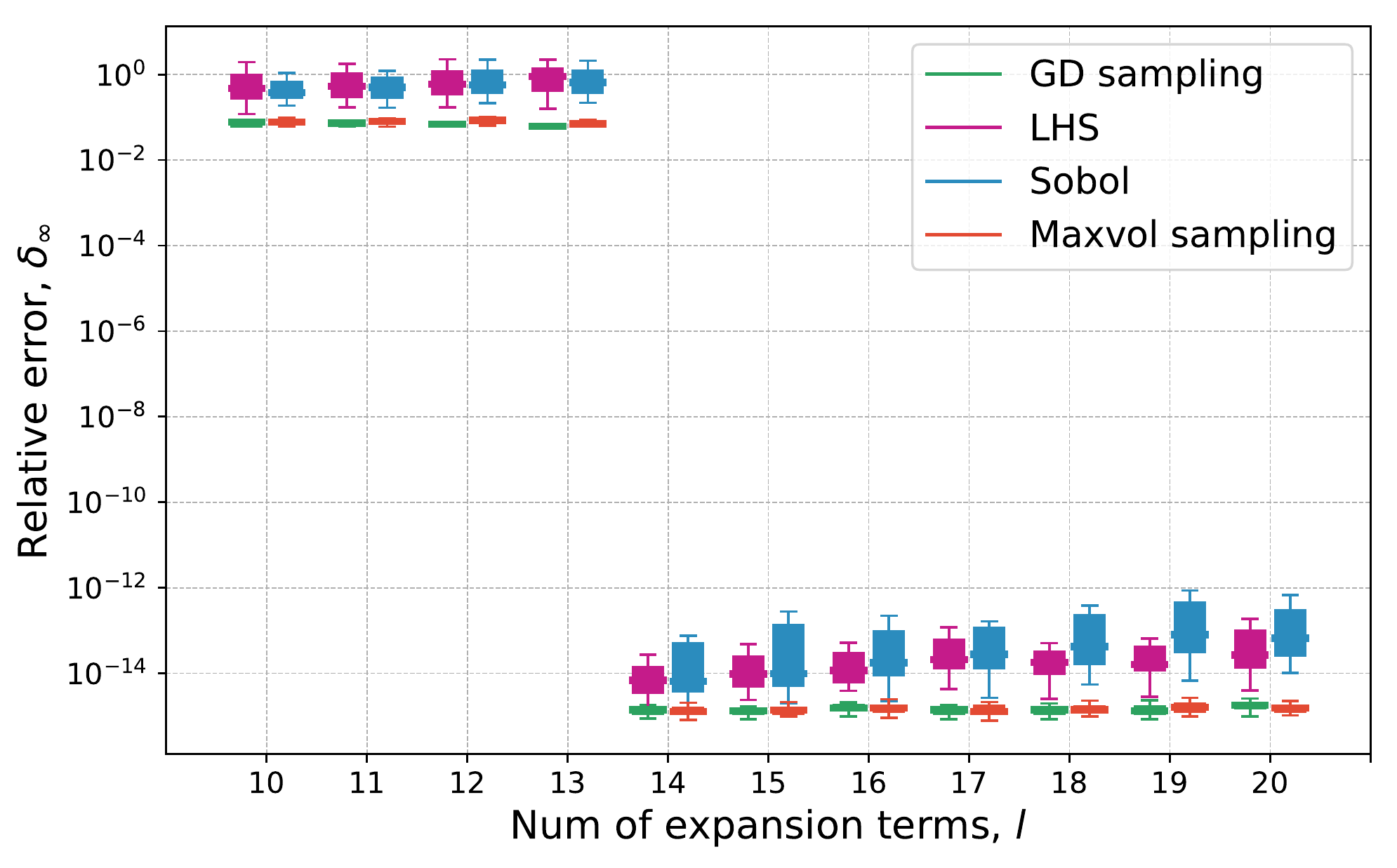}}\caption{}
    \end{subfigure}
    \hfill
    \begin{subfigure}[t]{0.49\linewidth}
        \centering{\includegraphics[width=1\linewidth]{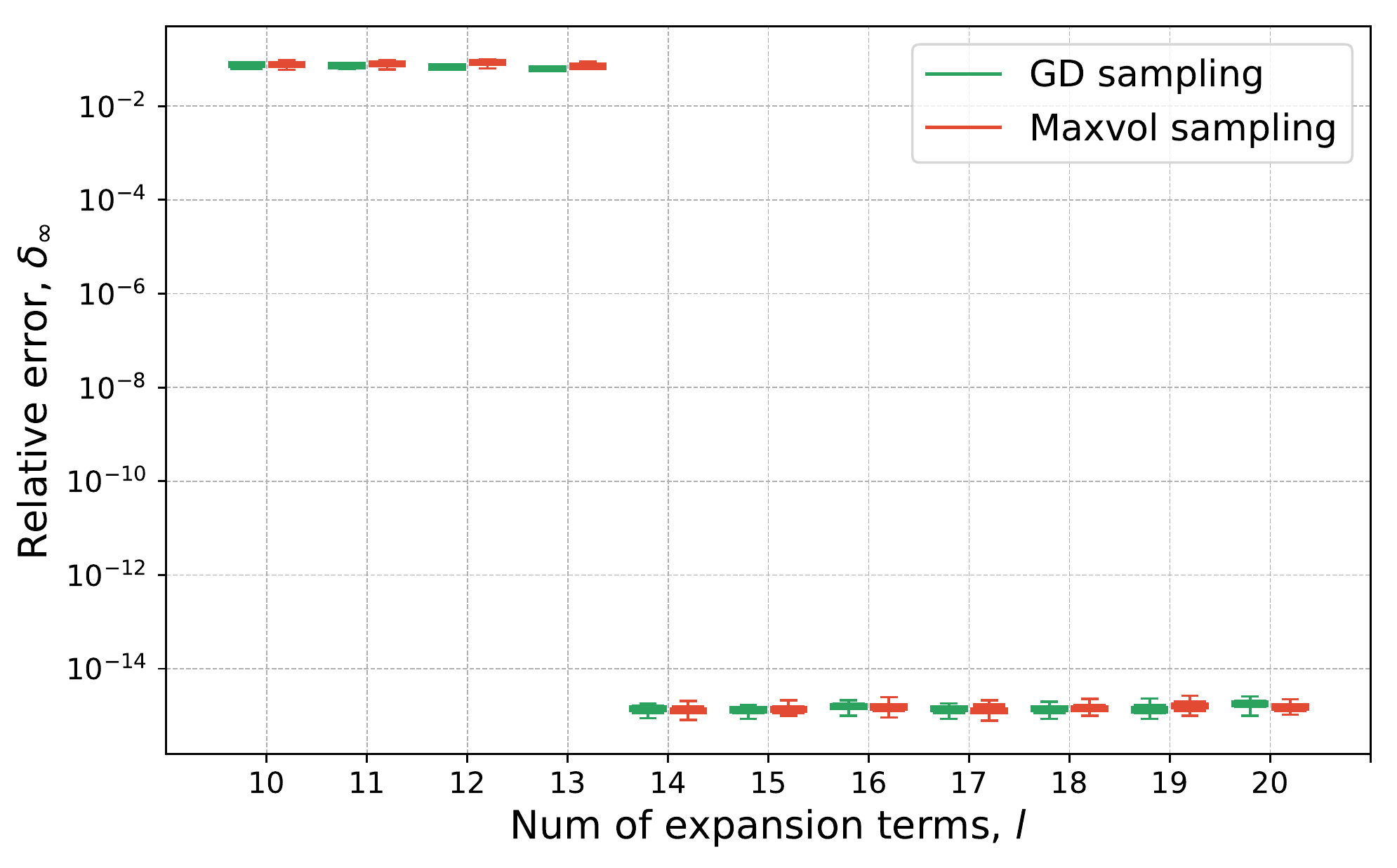}}\caption{}
    \end{subfigure}}
    \caption{\label{pic:4}\textbf{(a)}: Plot of Rosenbrock function on the domain of interest: $[-1,1]\times[-1,1]$. \textbf{(b)}: Test error for different number of terms (from $l = 10$ to $l = 20$) in polynomial expansion. The box-plots are obtained from 50 repetitions. \textbf{(c)}: The same as \textbf{(b)} but focused on two best sampling methods.}
\end{figure}

The performance of different sampling methods is compared in terms of infinity norm of the relative error on the test set for the varying size of experimental designs (\cref{pic:4}). Each analysis is repeated 50 times in order to estimate statistical uncertainty.

One can observe that ED based on LHS and Sobol' sequences show a poor performance compared to the GD sampling and Maxvol sampling. At the same time GD sampling and Maxvol sampling have the same performance. A drop in the accuracy from $l=13$ to $l=14$ for all the sampling methods is connected with the increase of the total degree of the polynomial expansion~$\widetilde f(\vec x)$, since the total degree of \cref{eq:rosenbrock} is 4, we get exact representation.

\subsubsection{Sincos function}
Now let us consider another two-dimensional function on a square $[-1,1]\times[-1,1]$ that we have denoted as \emph{sincos}: 
$$
f(x,y) = \sin\left(\frac{x^2}{2} - \frac{y^2}{4} + 3\right)\cdot\cos\left(2x+1-e^y\right).
$$

\begin{figure}[htb]
\centering{
    \begin{subfigure}[t]{0.37\textwidth}
        \centering{\includegraphics[width=1\textwidth]{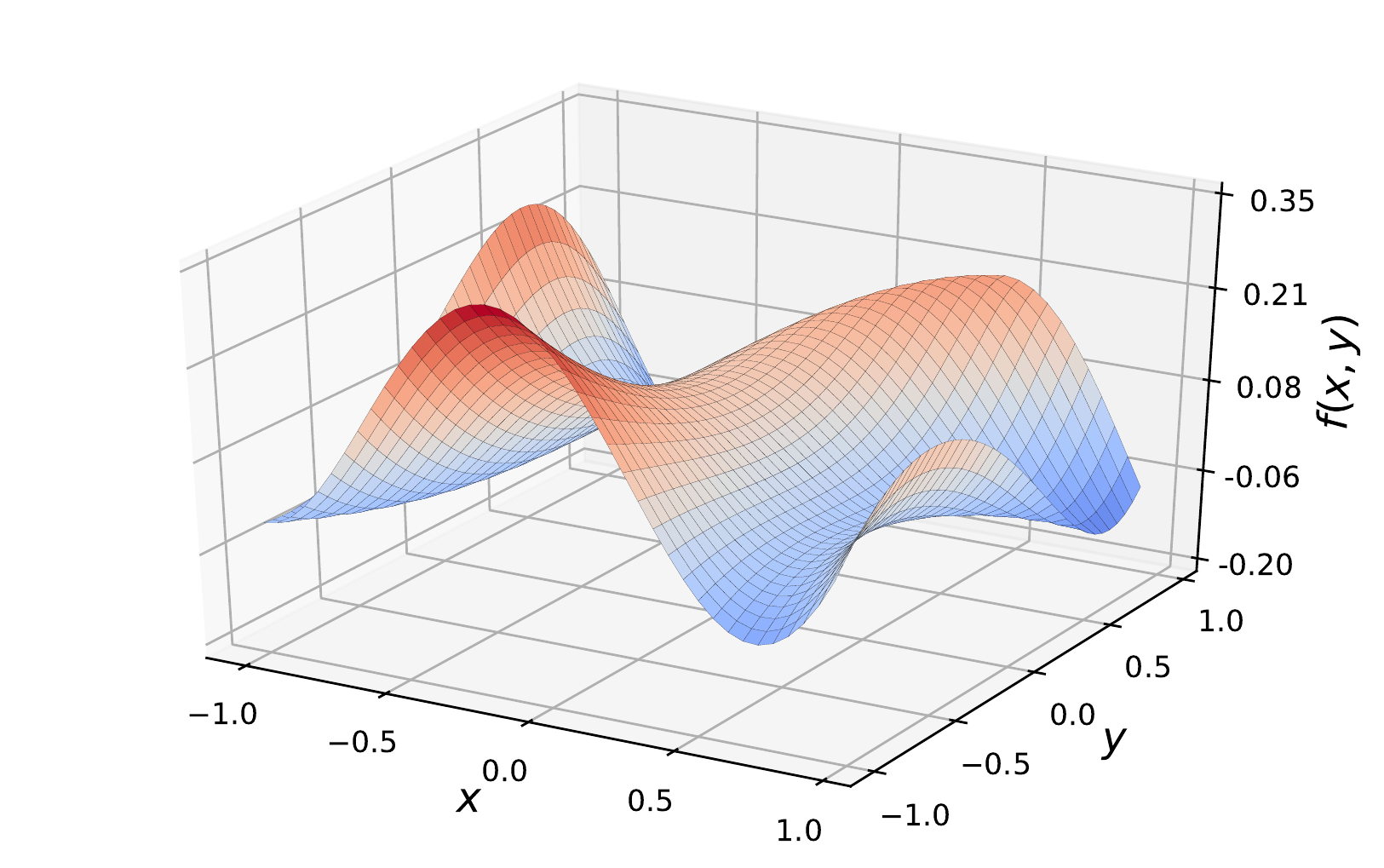}}\caption{Plot of sincos function.}
    \end{subfigure}
    \vfill
    \begin{subfigure}[t]{0.49\textwidth}
        \centering{\includegraphics[width=1\textwidth]{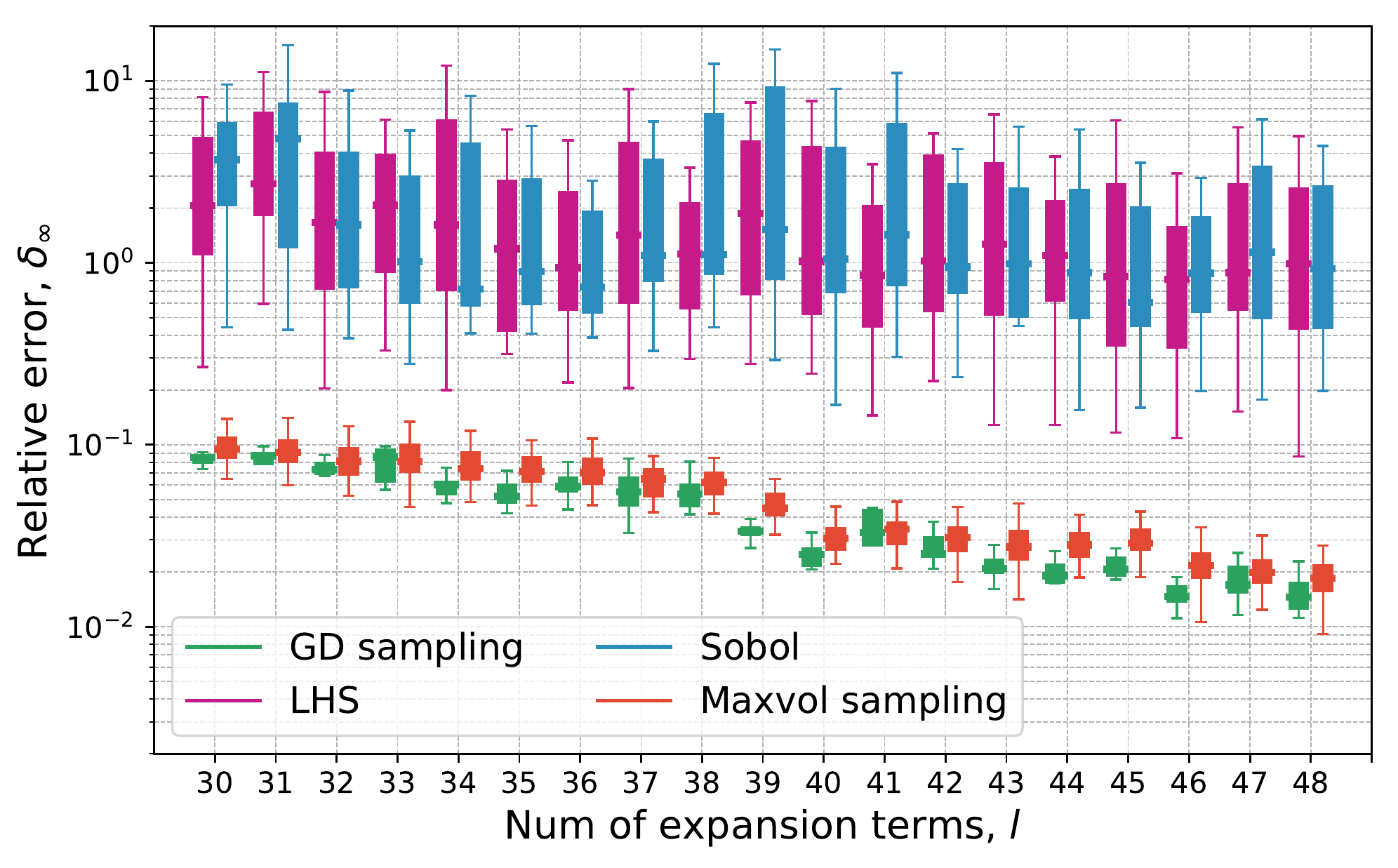}}\caption{}
    \end{subfigure}
    \hfill
    \begin{subfigure}[t]{0.49\textwidth}
        \centering{\includegraphics[width=1\textwidth]{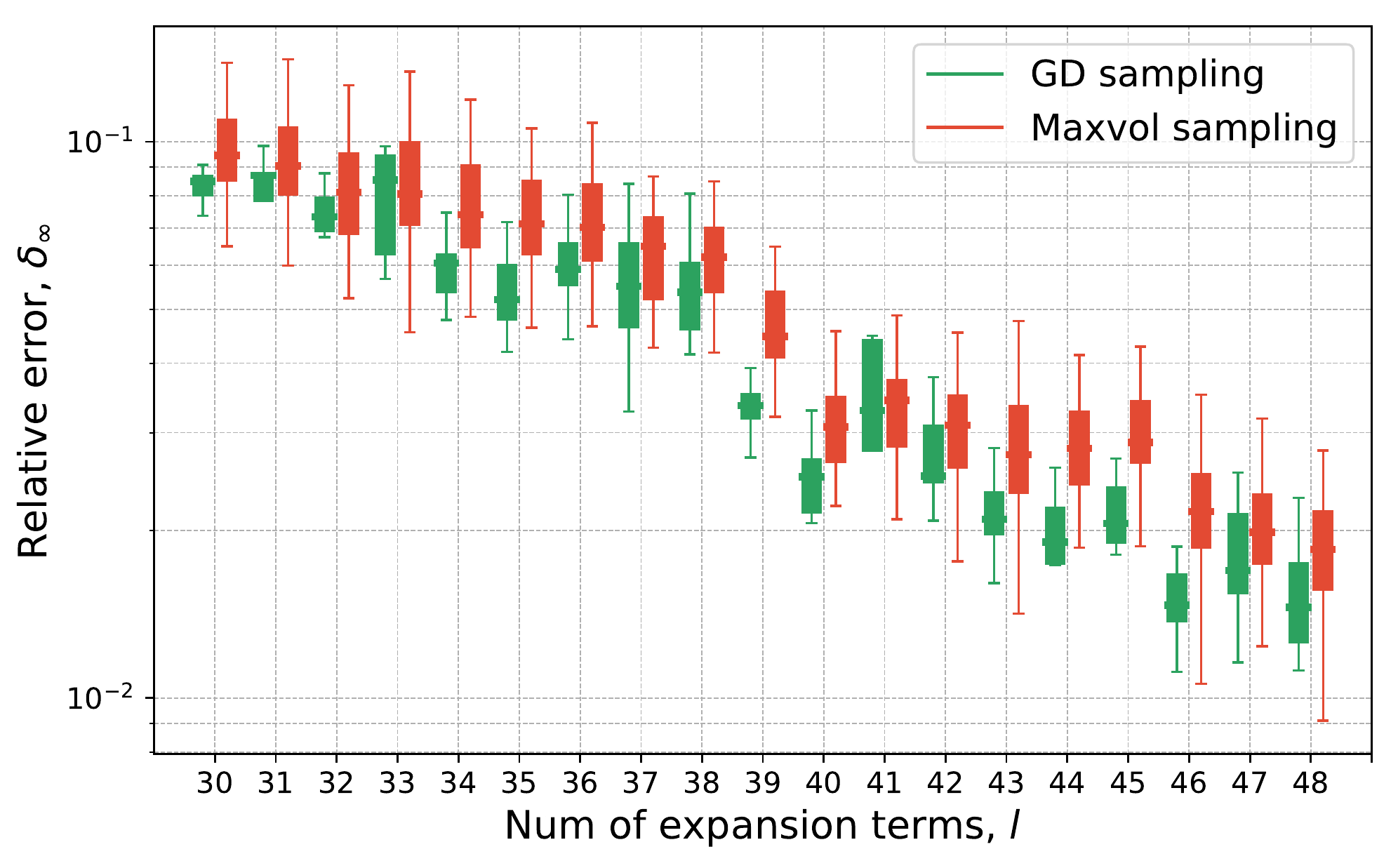}}\caption{}
    \end{subfigure}
    }\caption{\label{pic:5}\textbf{(a)}: Plot of sincos function on the domain of interest: $[-1,1]\times[-1,1]$. \textbf{(b)}: Evolution of box-plots of the test error for different number of terms in polynomial expansion. The box-plots are obtained from 50 repetitions. \textbf{(c)}: The same as \textbf{(b)} but focused on two best sampling methods.}
\end{figure}

As in the case of Rosenbrock function, for sincos function we compute the approximation error on $N_{\text{test}}=10^6$ test points for ED sizes in the range from $l = 30$ to $l = 48$ (\cref{pic:5}). On \cref{pic:5} the trend of decreasing approximation error with the increase of ED size for D-optimal sampling methods can be seen.     

\subsubsection{Gaussian function}
The final test on two-dimensional functions is a Gaussian function on the box domain $[-1,1]\times[-1,1]$:
$$f(x,y) = 2e^{-\frac{7}{2}\left(x^2+y^2\right)}.$$

\begin{figure}[htb]
\centering{
    \begin{subfigure}[t]{0.37\linewidth}
        \centering{\includegraphics[width=1\linewidth]{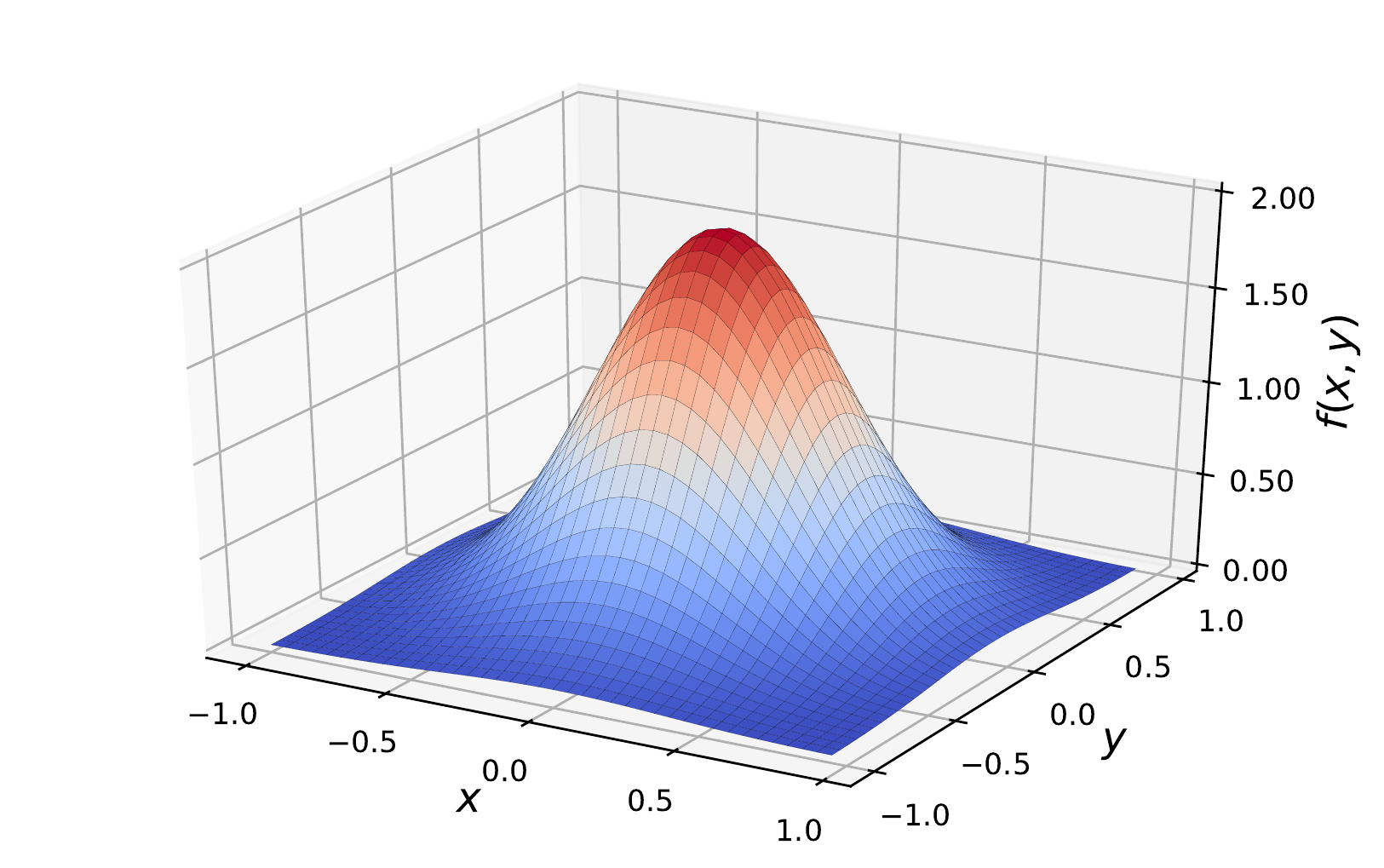}}\caption{Plot of Gaussian function.}
    \end{subfigure}
    \vfill
    \begin{subfigure}[t]{0.49\linewidth}
        \centering{\includegraphics[width=1\linewidth]{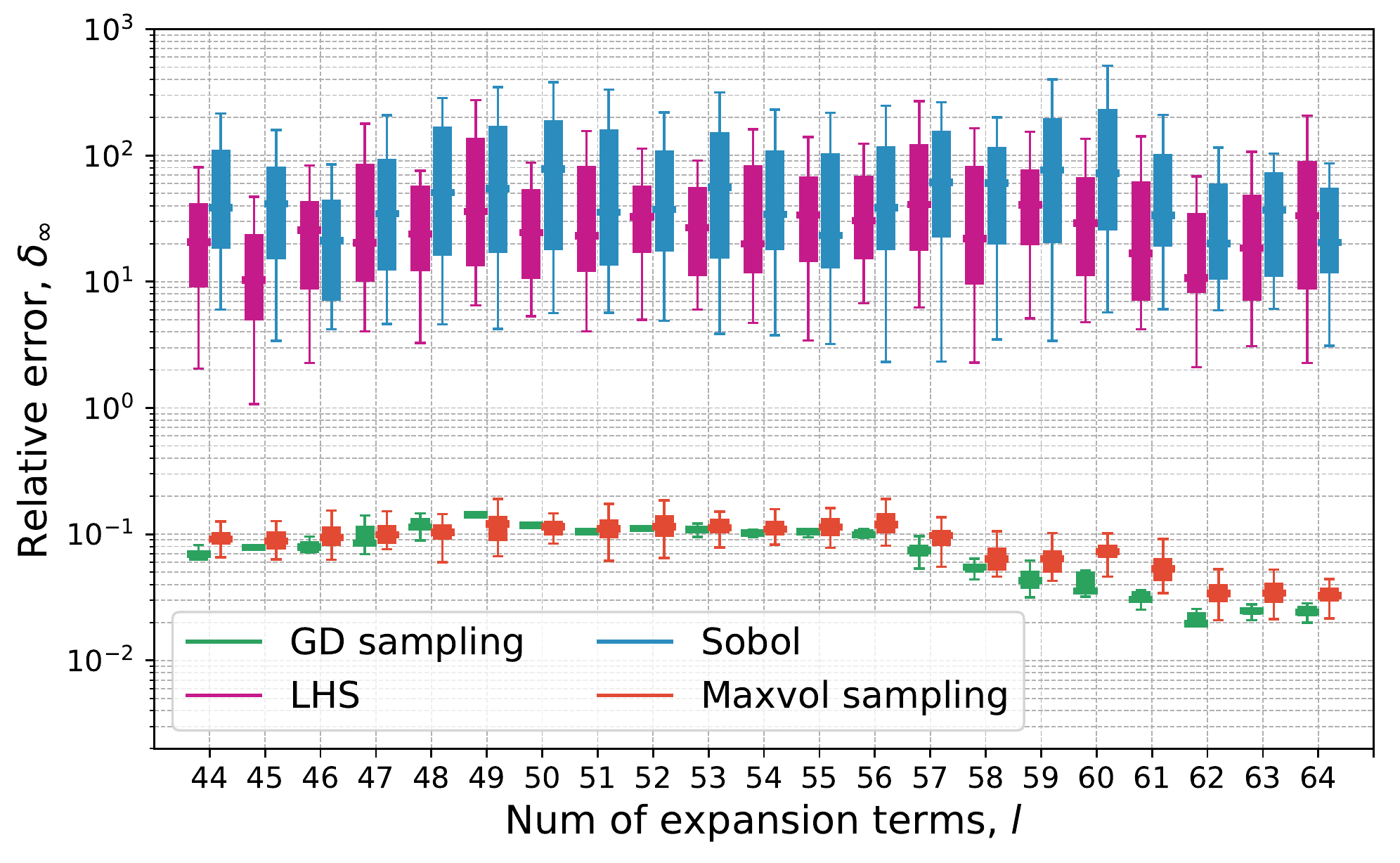}}\caption{}
    \end{subfigure}
    \hfill
    \begin{subfigure}[t]{0.49\linewidth}
        \centering{\includegraphics[width=1\linewidth]{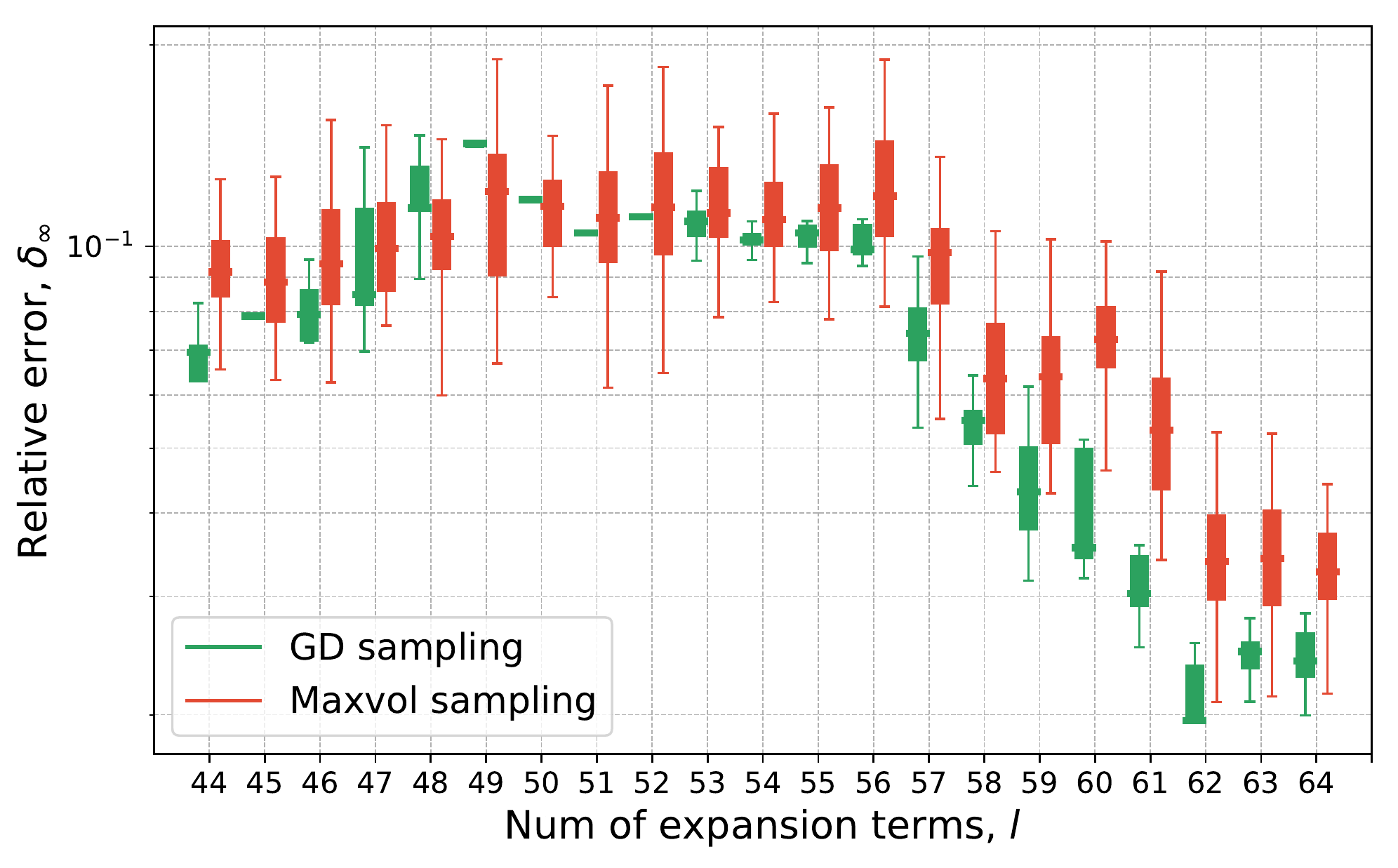}}\caption{}
    \end{subfigure}}\caption{\label{pic:6}\textbf{(a)}: Plot of Gaussian function on the domain of interest: $[-1,1]\times[-1,1]$. \textbf{(b)}: Evolution of box-plots of the test error for different number of terms in polynomial expansion. The box-plots are obtained from 50 repetitions. \textbf{(c)}: The same as \textbf{(b)} but focused on two best sampling methods.}
\end{figure}

On \cref{pic:6} one can see results similar to \cref{pic:5}: D-optimal sampling methods perform much better than LHS and Sobol' sequence sampling, and approximation error decreases with the increase in ED size. 

D-optimal design sampling methods consistently outperform other sampling methods. Moreover, such methods generally behave in a more stable way resulting in smaller variability between repetitions. Especially, this property becomes more important as the size of ED becomes larger.

\begin{figure}[ht]
\centering{
    \begin{subfigure}[t]{0.49\linewidth}
        \centering{\includegraphics[width=1\linewidth]{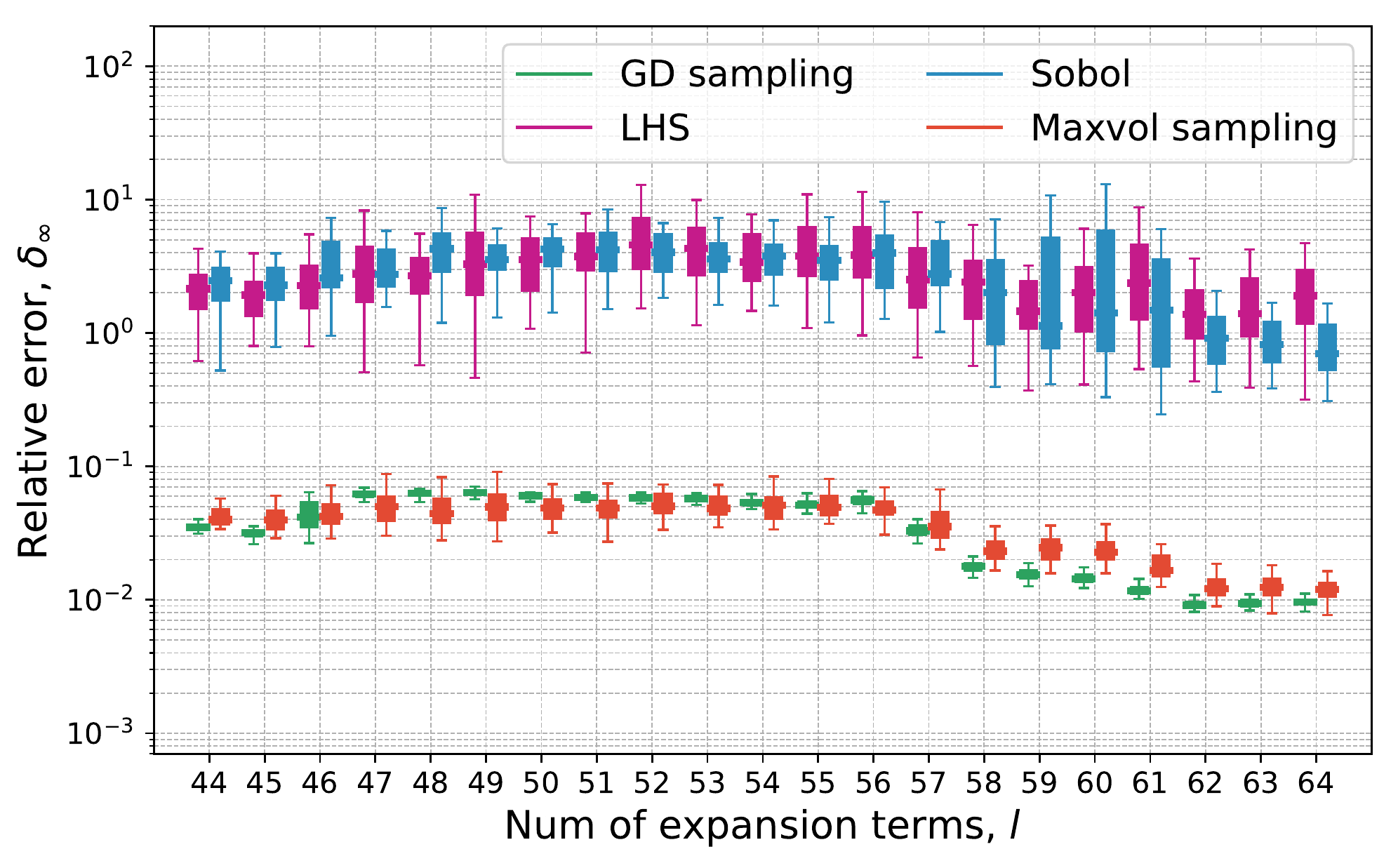}}\caption{\label{pic:7:a}}
    \end{subfigure}
    \hfill
    \begin{subfigure}[t]{0.49\linewidth}
        \centering{\includegraphics[width=1\linewidth]{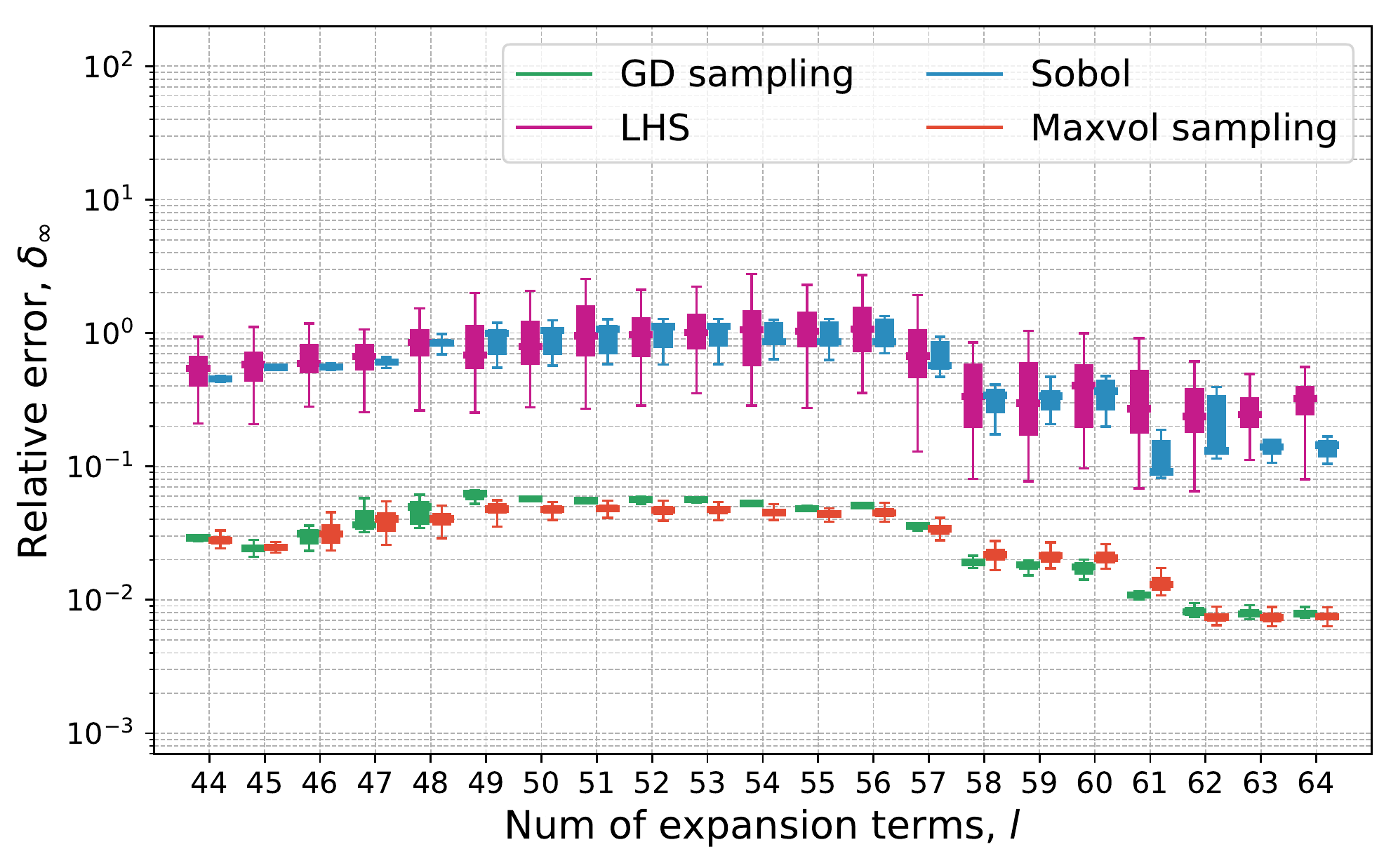}}\caption{\label{pic:7:b}}
    \end{subfigure}}
    \caption{\label{pic:7}Accuracy of approximation for different sizes of ED: \textbf{(a)} the number of sampled points is 1.1 times more than the number of terms (\textit{i.e.,} $n = \lceil1.1\cdot l\rceil)$, and \textbf{(b)} the number of sampled points is 2.5 times more than the number of terms in polynomial expansion (\textit{i.e.,} $n = \lceil2.5\cdot l\rceil)$.}
\end{figure}

Let us consider the case when more points are sampled than the number of terms in polynomial expansion (so-called \emph{oversampling}). 
So, we  complement the results for the approximation of Gaussian function (\cref{pic:6}) by considering two cases of oversampling: when the number of sampled points is 1.1 times more than the number of terms in polynomial expansion (\cref{pic:7:a}), and the case when we sample 2.5 times more points than the number of terms in corresponding expansion (\cref{pic:7:b}). As we can see on \cref{pic:7}, oversampling allows to improve the approximation (especially when compared to \cref{pic:6}) mainly for LHS and Sobol' sequence sampling. With the increase of the oversampling factor (from~1.1 to~2.5) Sobol' sequence sampling and GD sampling show more stable performance.  

Summarizing the results of the tests on two-dimensional functions, it can be stated that oversampling allows to significantly improve approximation accuracy for LHS and Sobol' sequence sampling, while practically has no effect on D-optimal sampling methods. 
It means that GD sampling or Maxvol sampling can be effectively used in the case of the tight budget on the number of runs of the complex model of interest. 
At the same time, it can be also noticed that oversampling provides more stable performance for GD sampling and Sobol' sequence sampling methods.    

\subsubsection{Piston simulation function}
In order to apply the proposed sampling method to high-dimensional surrogate modeling, we will consider a Piston simulation function~\cite{Kenett2013}. This function has a seven-dimensional input (see \cref{tbl:piston_descr}) and an one-dimensional output that effectively models the time (in seconds) that takes piston to complete one cycle within a cylinder. The cycle time is determined by a composition of functions: 

$$C(\vec{x}) = 2\pi\sqrt{\frac{M}{k+S^2\cdot\frac{P_0V_0}{T_0}\cdot\frac{T_a}{V^2}}},$$
$$\text{where}\ V = \frac{S}{2k}\left(\sqrt{A^2+4k\frac{P_0V_0}{T_0}T_a}-A\right),$$
$$A = P_0S+19.62M-\frac{kV_0}{S}.$$
\begin{table} [ht]
  \centering{
    \begin{tabular}{|p{1.3cm}||p{3.7cm}|p{2.4cm}|p{1cm}|}
     \hline
     Variable    & Name & Range & Units \\
     \hline
     $M$    & Piston weight           & [30, 60]         & $kg$\\
     $S$    & Piston surface area     & [0.005, 0.020]  & $m^2$\\
     $V_0$  & Initial gas volume      & [0.002, 0.010]  & $m^3$\\
     $k$    & Spring coefficient      & [1000, 5000]    & $N/m$\\
     $P_0$  & Atmospheric pressure    & [90000, 110000] & $N/m^2$\\
     $T_a$  & Ambient temperature     & [290, 296]       & $K$\\
     $T_0$  & Filling gas temperature & [340, 360]       & $K$\\
     \hline
    \end{tabular}\caption{\label{tbl:piston_descr}Description of Piston simulation variables.}
    }
\end{table} 

In the \cref{tbl:piston}, an approximation error (a median value over 30 runs) for two different sizes ($l=1750$ and $l=1850$) of ED and various number of sampling points~$n$ can be found.         

\begin{table}[htb]
  \centering
    \begin{tabular}{|p{3cm}||p{1.4cm}|p{1.4cm}|p{1.4cm}|p{1.4cm}|p{1.4cm}|}
     \hline
     \multicolumn{6}{|c|}{Number of expansion terms, $l = 1750$} \\
     \hline
     Number of samples    & $n = 1750$ & $n = 1760$ & $n = 1770$ & $n = 1780$ & $n = 1790$\\
     \hline
     $\delta_{\text{GD}}$      & \textbf{0.0702}     & \textbf{0.0593}     & \textbf{0.0577}     & \textbf{0.0486} & \textbf{0.0477}\\
     $\delta_{\text{Sobol}}$  & 0.5161     & 0.4676     & 0.3283    & 0.2656 & 0.2420\\
     $\delta_{\text{LHS}}$  & 0.3963     & 0.2494     & 0.3556     & 0.3481 & 0.2126 \\
     $\delta_{\text{Maxvol}}$  & 0.0928     & 0.0601     & 0.0526     & 0.0427 & 0.0400\\
     \hline
    \end{tabular}

\medskip 
  
    \begin{tabular}{|p{3cm}||p{1.4cm}|p{1.4cm}|p{1.4cm}|p{1.4cm}|p{1.4cm}|}
     \hline
     \multicolumn{6}{|c|}{Number of expansion terms, $l = 1850$} \\
     \hline
     Number of samples    & $n = 1850$ & $n = 1860$ & $n = 1870$ & $n = 1880$ & $n = 1890$\\
     \hline
     $\delta_{\text{GD}}$      & \textbf{0.0511}     & \textbf{0.0492}     & \textbf{0.0401}     & \textbf{0.0388} & \textbf{0.0352}\\
     $\delta_{\text{Sobol}}$  & 0.5925     & 0.5987     & 0.3521    & 0.2506 & 0.2517\\
     $\delta_{\text{LHS}}$  & 0.2357     & 0.6181     & 0.6009     & 0.3658 & 0.3169 \\
     $\delta_{\text{Maxvol}}$  & 0.0548     & 0.0489     & 0.0488     & 0.0462 & 0.0382\\
     \hline
    \end{tabular}
    \caption{\label{tbl:piston} Accuracy of the least-squares polynomial approximation of Piston simulation function.}
\end{table} 

From \cref{tbl:piston}, we can see that D-optimal sampling methods are the best. Also, it can be noted that oversampling allows to reduce the approximation error for the ED sampled with Sobol' sequence method.

Overall, it should be noted that in spite of similar performance of Maxvol sampling and GD sampling, the latter is more flexible as it is able to sample arbitrary points from the domain of interest whereas Maxvol sampling is limited with a discrete set of initial points.

\subsection{Lebesgue constant estimation}
In this Subsection we conduct a numerical estimation for the Lebesgue constant.

Let $\mathsf P(X)\colon C\rightarrow C$ be a projector on the span of the selected basis,
based on the experimental design matrix~$X$:
\begin{equation*}
(\mathsf P(X)f)(\vec x)=\widetilde f(\vec x),
\end{equation*}
where~$\widetilde f$ is  given by~\eqref{eq:1}
with coefficients defined in~\eqref{eq:LSM}.
By definition, Lebesgue constant~$\Lambda(X)$ is the $\infty$-norm of the operator~$\mathsf P(X)$:
\begin{equation*}
\Lambda(X)=\sup_{\normi|f|=1}\normi|\mathsf P(X)f|.
\end{equation*}
Let $f^*\in\Imm\P(X)$ be the best polynomial approximation of $f$
in $\infty$-norm,
so that $\normi|f^*-f|$ reaches minimum.
Then,
\begin{multline*}
\normi|f-\P(X)f|
\leq \normi|f-f^*| + \normi|f^*-\P(X)f|
=\normi|f-f^*| + \normi|\P(X)(f^*-f)|\\
\leq \normi|f-f^*| + \Lambda(X)\normi|f^*-f|
=(1+\Lambda(X))\normi|f-f^*|.
\end{multline*}
Thus, the Lebesgue constant can be utilized as an estimation of the approximation error obtained with our method in comparison with the best polynomial approximation of the same degree.

We will numerically estimate the Lebesgue constant using the same technique as for the estimation of the approximation error. Namely, we will take the maximum over the fixed set~$\mathcal D$ of points from the domain of interest:
\begin{equation*}
\Lambda(X)\approx 
\Lambda^\D(X)=\max_{\vec x\in\D}\sup_{\normi|f|=1}\abs|(\mathsf P(X)f)(\vec x)|.
\end{equation*}

To take the supremum over $f$, we apply the following trick.
Let us write the expansion~\eqref{eq:1} as a scalar product
of 
vector of basis functions $\Psi(\vec x)=\{\Psi_{\alpha_1}(\vec x),\,\Psi_{\alpha_2}(\vec x),\,\ldots\}$ and 
vector of coefficients~$c$ (see \Cref{eq:LSM}).
So, we obtain
\begin{equation*}
(\P(X)f)(\vec x)=\Psi(\vec x)\cdot  (A^TA)^{-1}A^T\mathcal Y.
\end{equation*}
Taking the supremum with respect to the function~$f$ is equivalent to taking the supremum over all vectors of $\mathcal Y$ such that $\norm|\mathcal Y|_1=1$, which in turn coincides with the first norm of the corresponding vector:
\begin{equation*}
\Lambda^\D(X)=\max_{\vec x\in\D}\norm|\Psi(\vec x)\cdot  (A^TA)^{-1}A^T|_1.
\end{equation*}

Estimation of the Lebesgue constant growth for ED obtained with particular sampling technique w.r.t. size of ED allows us to implicitly estimate accuracy of the least-squares polynomial approximation built using this ED. 

\begin{figure}[htb]
\centering{
    \begin{subfigure}[t]{0.49\linewidth}
        \centering{\includegraphics[width=\linewidth]{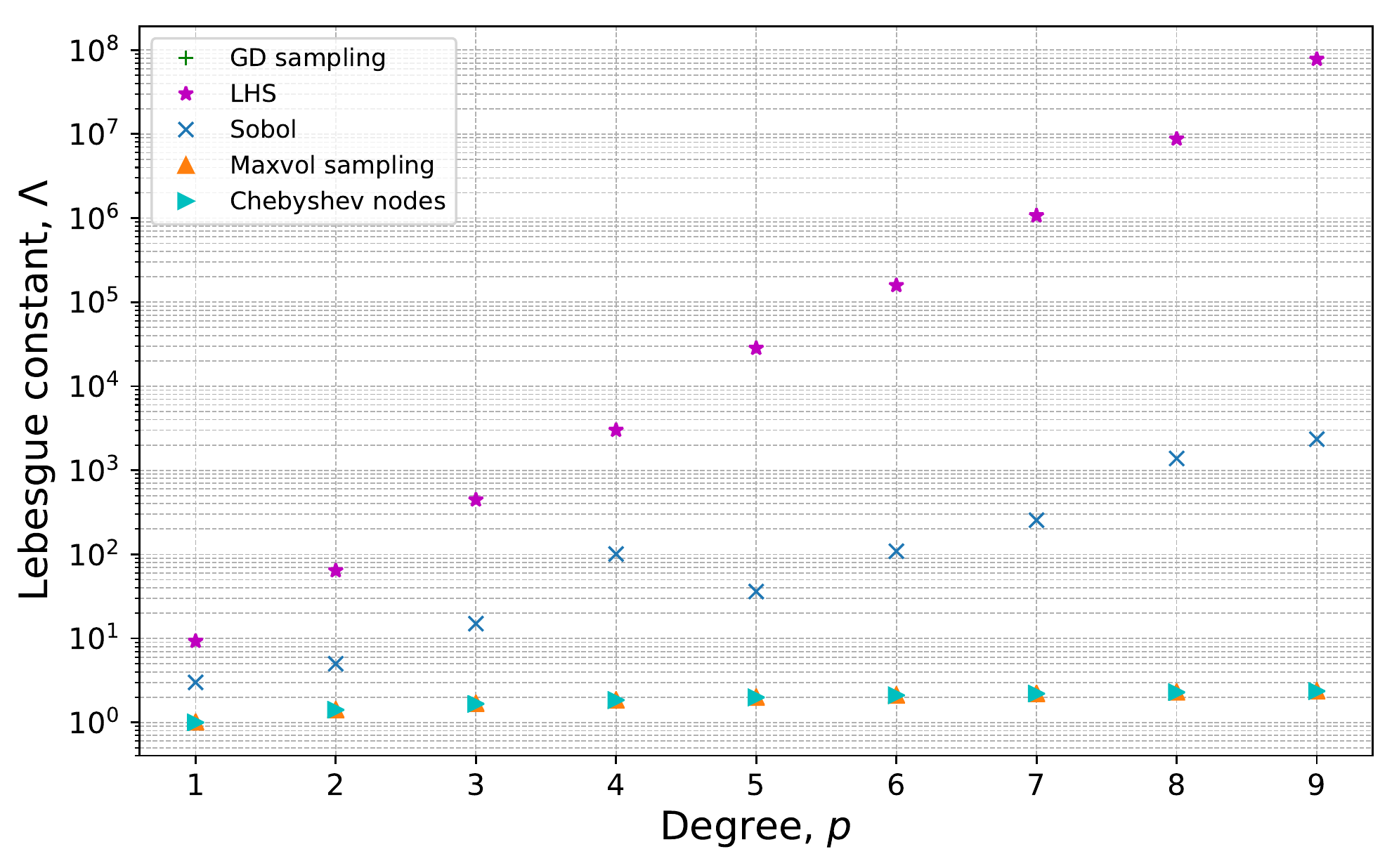}}
        \caption{\label{pic:8:a}}
    \end{subfigure}
    \hfill
    \begin{subfigure}[t]{0.49\linewidth}
        \centering{\includegraphics[width=\linewidth]{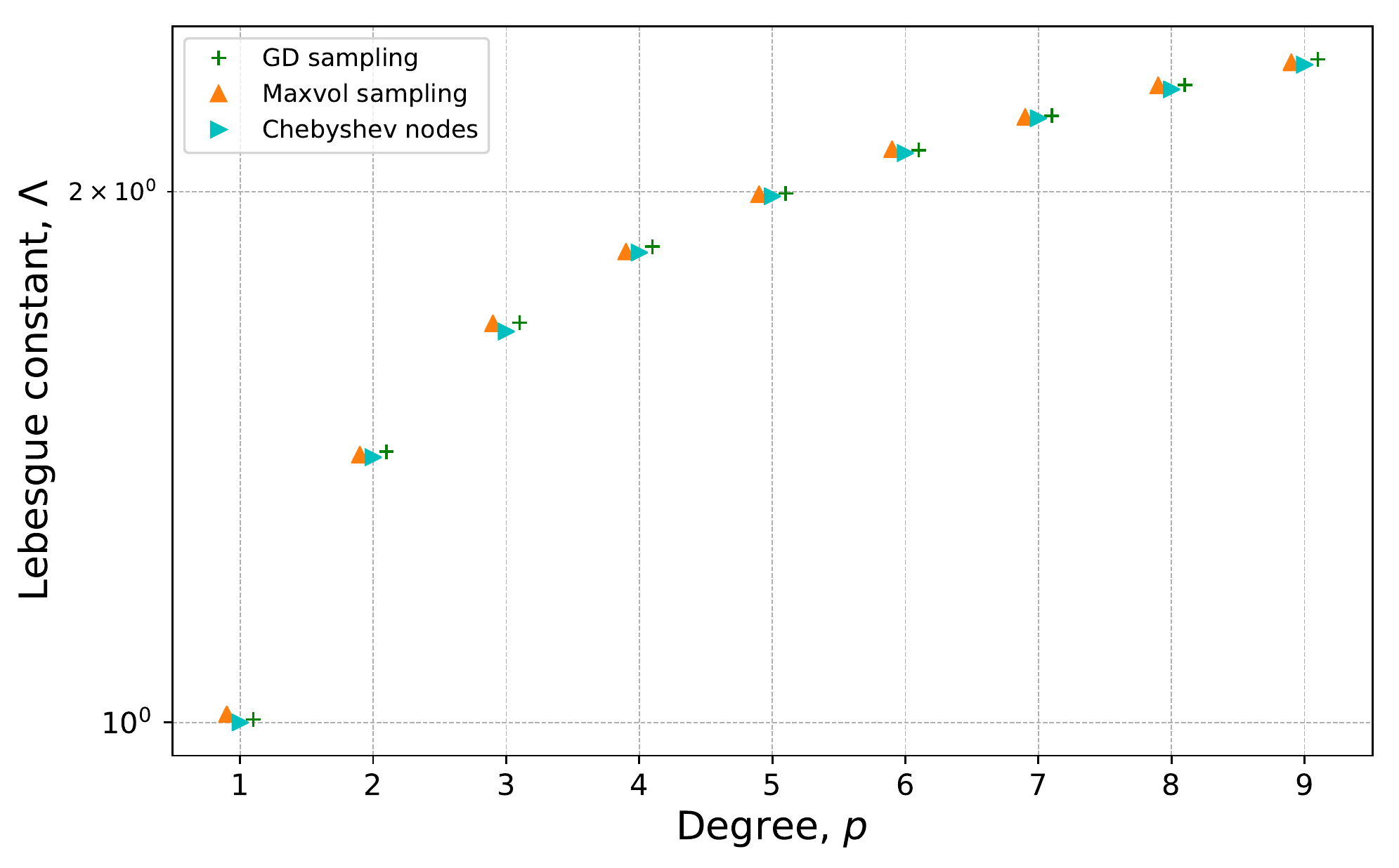}}
        \caption{\label{pic:8:b}}
    \end{subfigure}}
    \caption{\label{pic:8}\textbf{(a)}: Lebesgue constant estimation for the different ED sizes (\textit{i.e.} total polynomial degrees). \textbf{(b)}: The same data as \textbf{(a)} but exclusive LHS and Sobol' sequence.}
\end{figure}

At first, let us consider a one-dimensional case. It effectively means that the number of terms in polynomial expansion is equal to the total polynomial degree plus one: $l = p+1$. We use a test set of the size $N_{\text{test}}=10^6$ on the interval $[-1,1].$ We calculate the estimate of the Lebesgue constant $\Lambda_l$ for the points sampled by all the sampling techniques with respect to the number of such points in the range from $p=1$ to $p=9$. Also, we plot values of the Lebesgue constant for the Chebyshev roots as a quasi-optimal reference (\cref{pic:8}). 

As it can be seen from \cref{pic:8:a}, D-optimal sampling methods show much slower Lebesgue constant growth compared to LHS and Sobol' sequence. On the more detailed \cref{pic:8:b} one can find that the Lebesgue constant estimates for D-optimal sampling methods shows comparable growth with the Chebyshev nodes.

\begin{figure}[htb]
\centering{\centering{\includegraphics[width=0.9\linewidth]{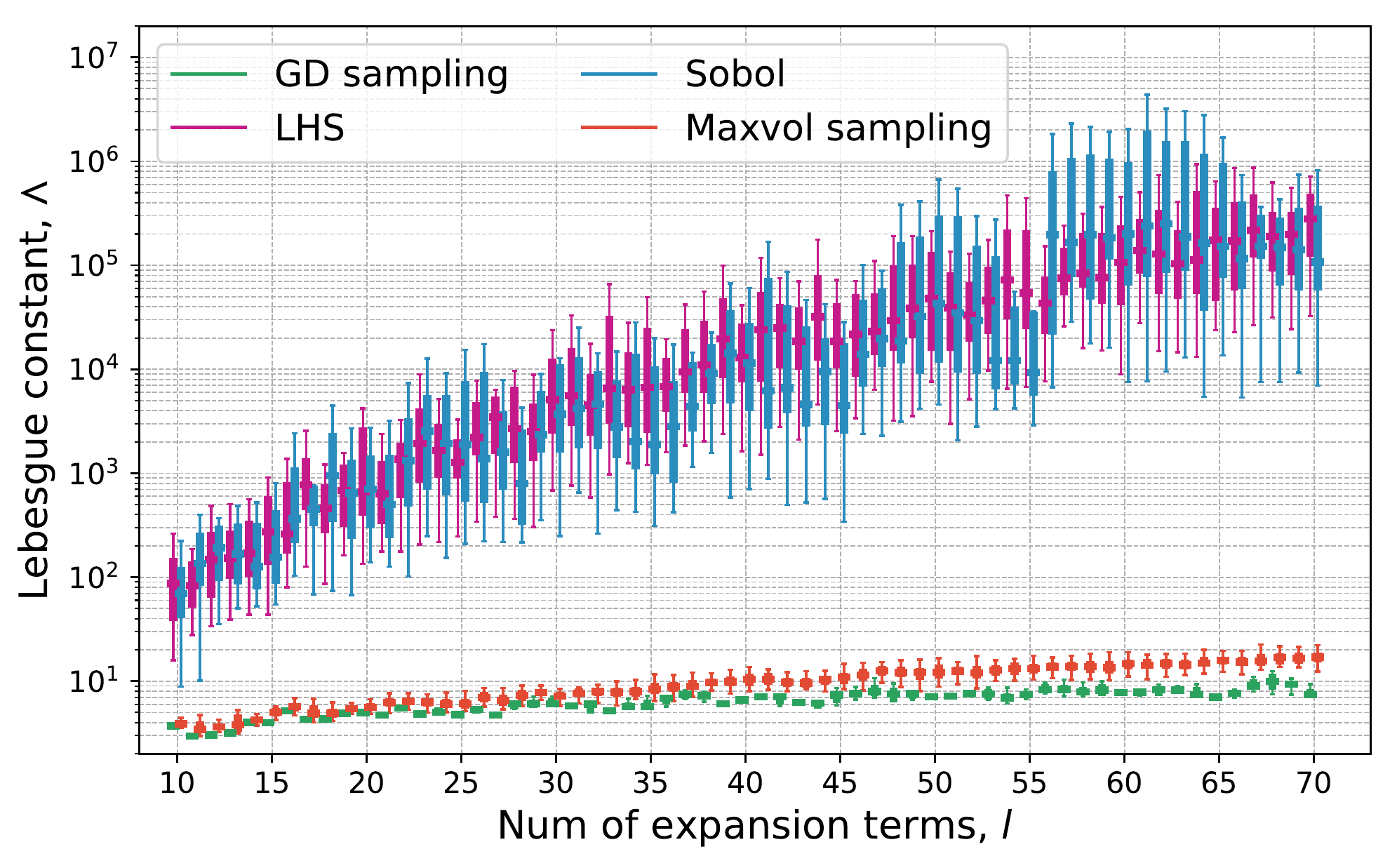}}}\caption{\label{pic:9}Lebesgue constant estimation for the ED size varying from 10 to 70.}
\end{figure}

We can also estimate the Lebesgue constant for the two-dimensional case. Results are shown on \cref{pic:9} where the size of ED is varying from 10 to 70 points. Since all sampling methods have a stochastic nature, for each ED size $l$ the main model was run 50 times and results were organized in a corresponding box-plot. 

From \cref{pic:9}, it can be seen that, as expected, in the two-dimensional case D-optimal sampling techniques perform much better than LHS and Sobol' sequence that is in consistence with the corresponding results on accuracy of the approximation (\cref{pic:4} -- \cref{pic:6}).
 
\subsection{Sampling from non-rectangular domain}
\begin{figure}[htb]
\centering{
    \begin{subfigure}[t]{0.32\linewidth}
        \centering{\includegraphics[width=1\linewidth]{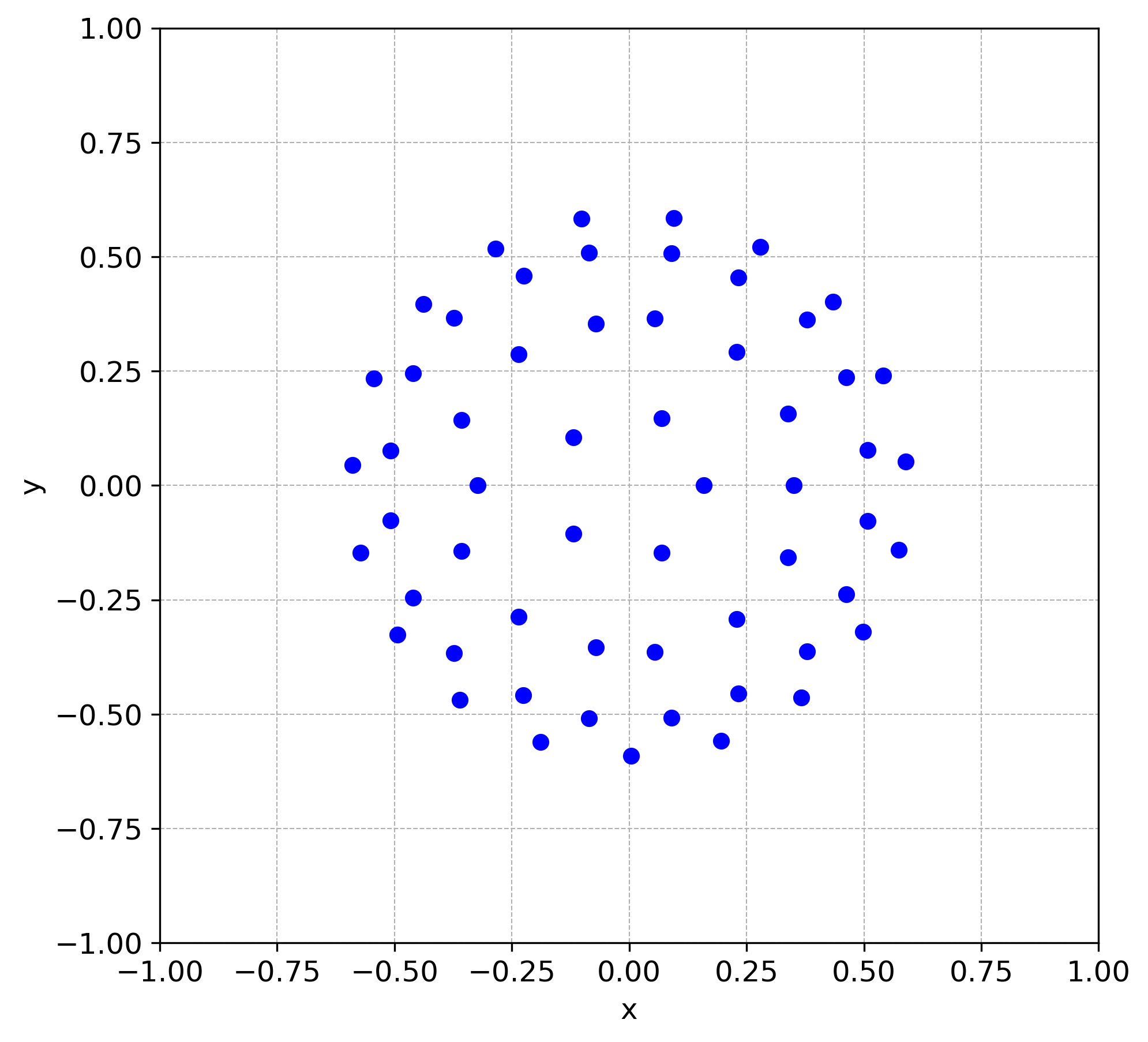}}\caption{Cirle-shape domain.}
    \end{subfigure}
        \hfill
    \begin{subfigure}[t]{0.32\linewidth}
        \centering{\includegraphics[width=1\linewidth]{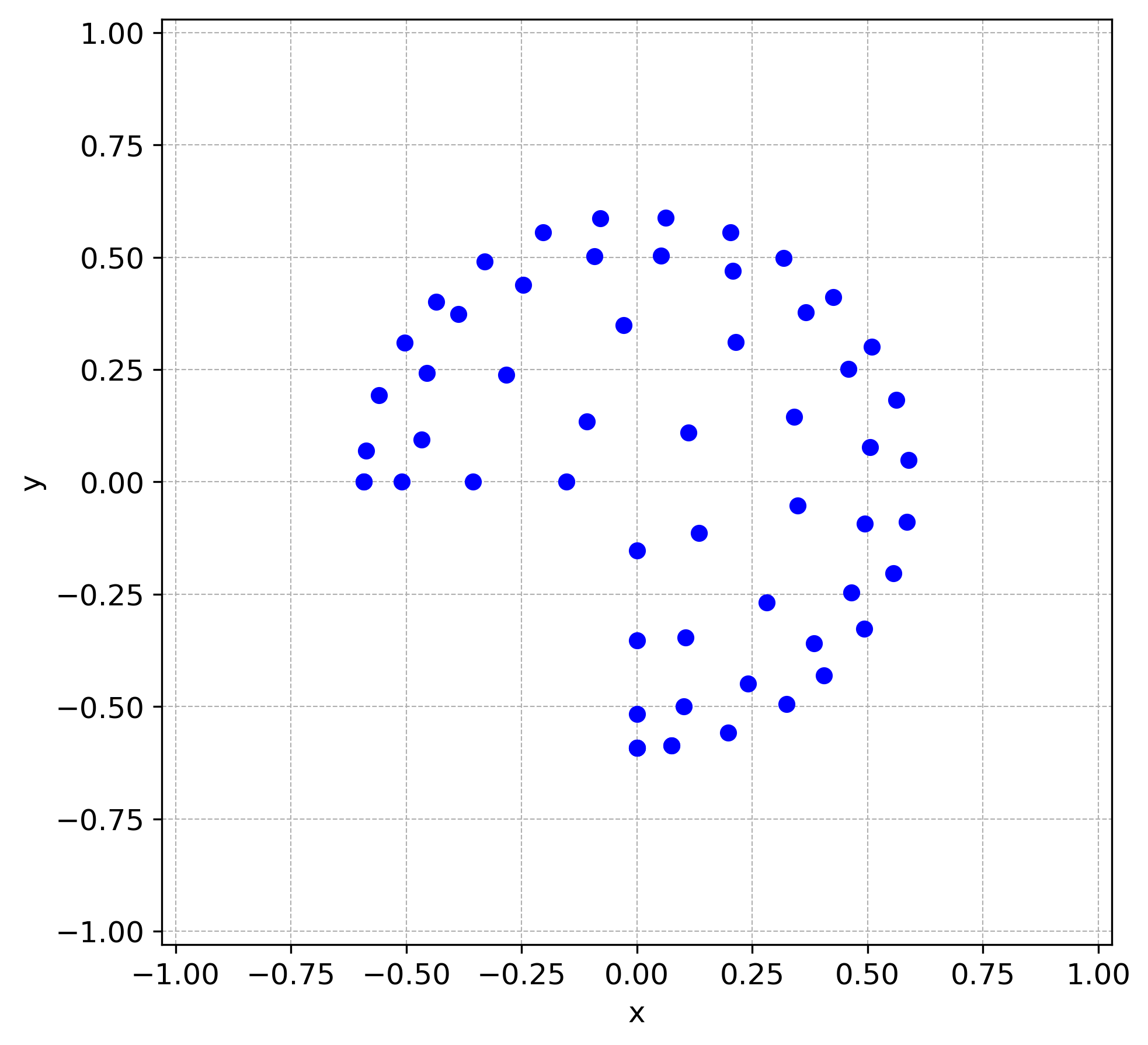}}\caption{Three quarters of a circle.}
    \end{subfigure}
    	\hfill
    \begin{subfigure}[t]{0.32\linewidth}
        \centering{\includegraphics[width=1\linewidth]{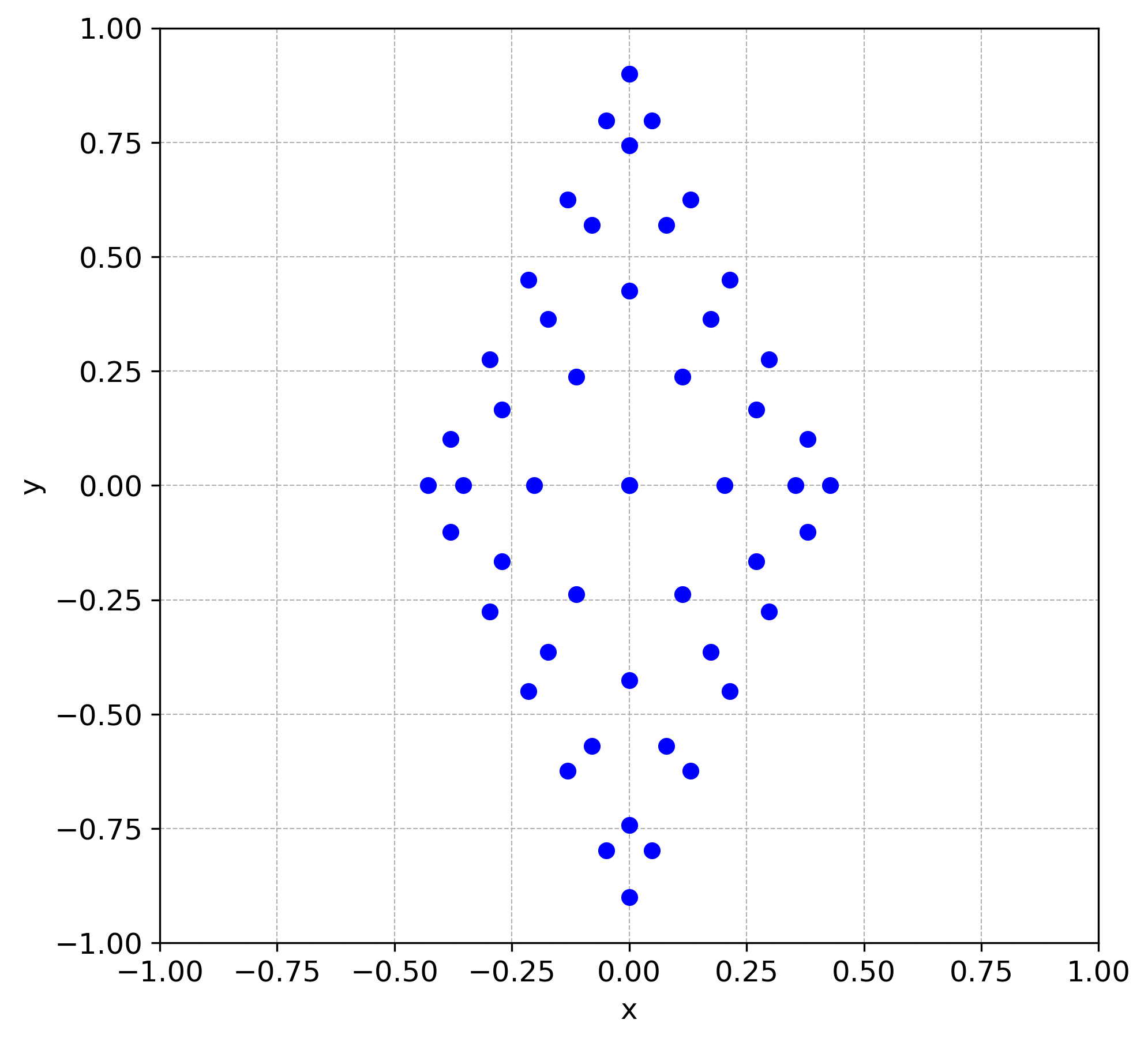}}\caption{Diamond-shape.}
    \end{subfigure}}\caption{\label{pic:10}EDs sampled from three non-rectangular two-dimensional domains for $n = 50$ and total degree $p = 5$.}
\end{figure}

In this subsection, we demonstrate an ability of the proposed method to sample points not only from rectangular domains but also from domains of arbitrary shape. The results showed on \cref{pic:10} were obtained with the use of \texttt{IPOP} optimizer~\cite{Wchter2005} for three two-dimensional domains with various shapes.  
As it can be seen from \cref{pic:10}, GD sampling shows quite a nice coverage of non-rectangular domains.

\section{Conclusions}
\label{sec:conclusion}

In this work, a new sampling method for finding a D-optimal experimental design is proposed. The developed method is based on the gradient descent algorithm and allows to minimize the log-det functional that determines the model matrix of corresponding experimental design.

The proposed sampling method was applied to the problem of the least-squares polynomial approximation of multivariate functions. Its efficiency was demonstrated by numerical comparison with the other sampling methods in the task of the multivariate function approximation. Numerical estimations on the Lebesgue constant growth were obtained and demonstrated quite a slow growth for the proposed sampling method.   

In the future, we plan to modify the proposed sampling method in order to make it adaptive and test it on the weighted least-squares polynomial approximation.

\appendix
\section{Proofs of \cref{lemma:one} and \cref{th:one}}
\label{sec:appendix}

At first, let us prove \cref{lemma:one}.
\begin{proof}
By the definition of $\widehat{W}(X):$
$$\dfrac{\partial\widehat{W}(X)}{\partial x_k} = -\dfrac{\partial\log\det B(X)}{\partial x_k}.$$
We use the chain rule of differentiation:
$$\dfrac{\partial\log\det B}{\partial x_k} =  \sum\limits_{i,j}\dfrac{\partial\log\det(B)}{\partial B_{ij}}\cdot\dfrac{\partial B_{ij}}{\partial x_k}.$$
First, we write down the first factor:
$$\dfrac{\partial\log\det(B)}{\partial B_{ij}} = \dfrac{\partial\det(B)}{\partial B_{ij}}\cdot\dfrac{1}{\det(B)} = \det(B)\cdot(B^{-1})_{ji}\cdot\dfrac{1}{\det(B)} = (B^{-1})_{ji}.$$
Let us write down the second factor:
\begin{multline*}
\dfrac{\partial B_{ij}}{\partial x_k} = \dfrac{\partial\left(\sum\nolimits_{m=1}^{n}A_{mi}A_{mj}\right)}{\partial x_k} = \dfrac{\partial\left(\sum\nolimits_{m=1}^{n}\psi_{i}(x_m) \psi_{j}(x_m)\right)}{\partial x_k} = \dfrac{\partial\left(\psi_{i}(x_k) \psi_{j}(x_k)\right)}{\partial x_k} \\
= \dfrac{\partial \psi_i(x_k)}{\partial x_k}\cdot\psi_j(x_k) + \psi_i(x_k)\cdot\dfrac{\partial \psi_j(x_k)}{\partial x_k}.    
\end{multline*}
Thus, we finally obtain the indicated equality:
$$G_k = \dfrac{\partial\widehat{W}(X)}{\partial x_k} =  -\sum\limits_{ij}(B^{-1}(X))_{ji}\cdot\left[\dfrac{\partial \psi_i(x_k)}{\partial x_k}\cdot\psi_j(x_k) + \psi_i(x_k)\cdot\dfrac{\partial \psi_j(x_k)}{\partial x_k}\right].$$
\end{proof}

Now, let us prove \cref{th:one}.
\begin{proof}
Using \cref{lemma:one}, we get the following expression:
$$G_{kl} = \dfrac{\partial\widehat{W}(X)}{\partial x_k^{(l)}} = -\sum\limits_{i,j}(B^{-1}(X))_{ji}\cdot\dfrac{\partial B_{ij}(X)}{\partial x_k^{(l)}}.$$
We will separately write down the second factor.
Note, that matrix element~$A_{mi}$
depends only on the point~$\vec x_m$:
$A_{mi}=A_{mi}(\vec x_m)$.
The derivatives of this element vanish along the components of the remaining points:
\begin{equation*}
\frac{\partial B_{ij}(X)}{\partial x_k^{(l)}} =
\frac{\partial\left(\sum_{m=1}^nA_{mi}(\vec x_m)A_{mj}(\vec x_m)\right)}{\partial x_k^{(l)}}=
\frac{\partial\left(A_{ki}(\vec x_k)A_{kj}(\vec x_k)\right)}{\partial x_k^{(l)}}.
\end{equation*}
Now, we use the explicit form of \Cref{eq:vandermonde}: 
\begin{multline*}
\dfrac{\partial B_{ij}(X)}{\partial x_k^{(l)}} =
\frac{\partial\left(A_{ki}(\vec x_k)A_{kj}(\vec x_k)\right)}{\partial x_k^{(l)}}
=\frac{\partial}{\partial x_k^{(l)}}
\left(\prod\nolimits_{q=1}^{d}\psi_{\vec\alpha_i^{(q)}}(x_k^{(q)})\cdot\prod\nolimits_{s=1}^{d}\psi_{\vec\alpha_j^{(s)}}(x_k^{(s)})\right)\\
=\frac{\partial}{\partial x_k^{(l)}}\left(\prod\nolimits_{q=1}^{d}\psi_{\vec\alpha_i^{(q)}}({x}_k^{(q)})\cdot\psi_{\vec\alpha_j^{(q)}}({x}_k^{(q)})\right)\\
= \dfrac{\partial\left(\psi_{\vec\alpha_i^{(k)}}({x}_k^{(l)})\cdot\psi_{\vec\alpha_j^{(k)}}({x}_k^{(l)})\right)}{\partial {x}_k^{(l)}}\cdot\prod\limits_{\substack{q=0\\(q\neq{l})}}^{d-1}\psi_{\vec\alpha_i^{(q)}}({x}_k^{(q)})\cdot\psi_{\vec\alpha_j^{(q)}}({x}_k^{(q)})\\
= \left[\dfrac{\partial\psi_{\vec\alpha_i^{(k)}}({x}_k^{(l)})}{\partial {x}_k^{(l)}}\cdot\psi_{\vec\alpha_j^{(k)}}({x}_k^{(l)}) + \psi_{\vec\alpha_i^{(k)}}({x}_k^{(l)})\cdot\dfrac{\partial \psi_{\vec\alpha_j^{(k)}}({x}_k^{(l)})}{\partial {x}_k^{(l)}}\right]\\
\times\prod\limits_{\substack{q=0\\(q\neq{l})}}^{d-1}\psi_{\vec\alpha_i^{(q)}}({x}_k^{(q)})\cdot\psi_{\vec\alpha_j^{(q)}}({x}_k^{(q)}).
\end{multline*}
\end{proof}

\bibliographystyle{siamplain}
\bibliography{references}

\end{document}